\documentclass{amsart}
\usepackage{amscd,amssymb,graphics}

\newtheorem{theorem}{Theorem}[section]
\newtheorem{lemma}[theorem]{Lemma}
\newtheorem{proposition}[theorem]{Proposition}
\newtheorem{corollary}[theorem]{Corollary}

\theoremstyle{definition}
\newtheorem{definition}[theorem]{Definition}

\theoremstyle{remark}
\newtheorem{remark}[theorem]{Remark}

\numberwithin{equation}{section}

\newcommand{\NN}{\mathbb{N}}
\newcommand{\ZZ}{\mathbb{Z}}

\newcommand{\RR}{\mathbb{R}}
\newcommand{\CC}{\mathbb{C}}

\newcommand{\PP}{\mathbb{P}}

\newcommand {\shExt}  {\mathcal{E} \!\text{\textit{xt}}}

\newcommand {\shI}   {\mathcal{I}}

\newcommand {\shL}   {\mathcal{L}}

\newcommand {\shT}   {\mathcal{T}}

\newcommand {\shX}   {\mathcal{X}}

\newcommand {\foB}   {\mathfrak{B}}
\newcommand {\foC}   {\mathfrak{C}}

\newcommand {\foK}   {\mathfrak{K}}
\newcommand {\foS}   {\mathfrak{S}}
\newcommand {\foM}   {\mathfrak{M}}

\newcommand {\foX}   {\mathfrak{X}}
\newcommand {\foY}   {\mathfrak{Y}}
\newcommand {\foZ}   {\mathfrak{Z}}

\newcommand {\Aut}   {\operatorname{Aut}}

\newcommand {\Ext}   {\operatorname{Ext}}

\newcommand {\GL}   {\operatorname{GL}}

\newcommand {\Hom}   {\operatorname{Hom}}
\newcommand {\id}   {\operatorname{id}}

\renewcommand {\ker }  {\operatorname{kern}}

\newcommand {\lra}   {\longrightarrow}

\newcommand {\NS}   {\operatorname{NS}}

\renewcommand{\O}    {\mathcal{O}}

\newcommand {\Pic}   {\operatorname{Pic}}
\newcommand {\PGL}   {\operatorname{PGL}}

\newcommand {\pr}   {\operatorname{pr}}

\newcommand {\quadand} {\quad\text{and}\quad}

\newcommand {\ra}   {\rightarrow}

\newcommand {\Sing}  {\operatorname{Sing}}
\newcommand {\SL}   {\operatorname{SL}}
\newcommand {\Spec}  {\operatorname{Spec}}

\begin{document}

\title[Degenerations of Kodaira surfaces] {Irreducible degenerations
  of primary Kodaira surfaces}

\author[Stefan Schroeer]{Stefan Schr\"oer} \address{Mathematische
  Fakult\"at, Ruhr-Universit\"at, 44780 Bochum, Germany}
\email{s.schroeer@ruhr-uni-bochum.de}

\author{Bernd Siebert} \address{Mathematische Fakult\"at,
  Ruhr-Universit\"at, 44780 Bochum, Germany}
\email{bernd.siebert@ruhr-uni-bochum.de}



\dedicatory{Hans Grauert zum 70.\ Geburtstag gewidmet}
\begin{abstract}
  We classify irreducible $d$-semistable degenerations
  of primary Kodaira surfaces. As an application we construct a canonical
  partial completion for the moduli space of primary Kodaira surfaces.
\end{abstract}

\maketitle

\section*{Introduction}

A smooth compact complex surface with trivial canonical bundle is a
K3 surface, a 2-dimensional complex torus, or a primary Kodaira
surface. Normal crossing degenerations of such surfaces have
attracted much attention. For example, Kulikov analyzed projective
degenerations of K3 surfaces \cite{Kulikov 1977}. His results were
generalized by Persson and Pinkham to degenerations of surfaces with
trivial canonical bundle whose central fiber has algebraic components
\cite{Persson; Pinkham 1981}. Conversely, Friedman characterized the
singular K3 surfaces from Kulikov's list that deform to smooth K3
surfaces by his notion of $d$-semistability \cite{Friedman 1983}.

For non-K\"ahler degenerations no general results seem to be known.
We chose primary Kodaira surfaces as our object of study because they
lack some complications that the more interesting classes of tori and
K3 surfaces have. On the other hand, Kodaira surfaces have a lot in
common with tori and, via the Kummer construction, with K3 surfaces.
The phenomena one sees in degenerations of Kodaira surfaces should
therefore also be observable in the other classes of surfaces with
trivial canonical bundle.

By a result of Borcea the moduli space $\foK$ of smooth Kodaira
surfaces is isomorphic to a countable union of copies of
$\CC\times\Delta^*$ \cite{Borcea 1984}. The parameters correspond to
the $j$-invariant of the elliptic base and a refined $j$-invariant
of the fiber. As $\Delta^*$ can not be completed to a closed Riemann
surface by adding finitely many points, it is clear that by studying
degenerations we can at most hope for a partial completion of $\foK$.
The explicit form of the moduli space also suggests the existence of
families leading to a degeneration of the base or the fiber to a
nodal elliptic curve. We will see that this is indeed the case. This
leads to a partial completion $\CC\times\Delta^* \subset \PP^1\times
\Delta\setminus (\infty,0)$ of each component of $\foK$. However, as
elliptic curves can degenerate to any $k$-cycle of rational curves,
this is not the full story. Rather we obtain a whole hierarchy of
such completions, that are linked by non-Hausdorff phenomena at the
boundary. Moreover, we will see that at the most interesting point
$(\infty,0)$ the picture becomes complicated. We were not able to
fully clarify what happens there. If one restricts to normal crossing
surfaces one should certainly blow up this point. On the other hand,
we will also make the presumably not so surprising observation that,
as in the case of abelian varieties \cite{Alexeev 1999}, it does not
suffice to restrict to normal crossing varieties. We did however find
families of generalized Kodaira surfaces mapping properly to
$\Delta\times \PP^1$, whose singularities are at most products of
normal crossing singularities.

The bulk of the paper is concerned with a classification of
irreducible, $d$-semi\-stable, locally normal crossing surfaces $X$
with $K_X=0$. Three different types occur. Our main result is a
description of the resulting completion $\foK\subset\overline{\foK}$,
derived by deformation theory. The three different types correspond
to three different parts in the boundary $\foB= \overline{\foK}
\setminus\foK$. Each part is a countable union of copies of
$\Delta^*$ or $\CC$. Locally along the boundary divisor the
completion $\foK\subset\overline{\foK}$ looks like the blowing-up of
$(\infty,0)\in\PP^1\times\Delta$, with two points on the exceptional
divisor removed.

Some examples for degenerations of Kodaira surfaces have previously
been given by Friedman and Shepherd-Barron in \cite{Friedman;
Shepherd-Barron 1983}.

\medskip This article is divided into six sections. The first section
contains general facts about smooth Kodaira surfaces. In the second
section we describe their potential degenerations and show that three
types are possible. Sections 3--5 contain an analysis of each type. In
the final section we assemble our results in terms of moduli spaces,
complemented by some examples of smoothable surfaces with
singularities that are products of normal crossings.

We thank the referee for suggestions concerning Theorem~3.10 and the
interpretation of the completed moduli space after
Proposition~\ref{boundary components}.

\medskip This paper is dedicated to Hans Grauert on the occasion of
his 70th birthday. The second author wants to take this opportunity
to express his gratitude for the support and mathematical stimulus he
received from him as one of his last students. It was a great
pleasure to learn from him.

\section{Smooth Kodaira surfaces}

In this section we collect some facts on primary Kodaira surfaces. Suppose
$B,E$ are two elliptic curves, and endow $E$ with a group structure. Let
$f:X\ra B$ be a holomorphic principal $E$-bundle. The canonical bundle
formula gives $K_X=0$, so the Kodaira dimension is $\kappa(X)=0$. As a
topological space, $X$ is the product of the $1$-sphere with a
$1$-sphere-bundle $g:M\ra B$. Let $e(g)\in H^2(B,\ZZ)$ be its Euler class.
The homological Gysin sequence
$$
H_2(B,\ZZ) \stackrel{e(g)}{\lra} H_0(B,\ZZ) \lra H_1(M,\ZZ)
\lra H_1(B ,\ZZ) \lra 0
$$
and the K\"unneth formula yield $H_1(X,\ZZ)=\ZZ^3\oplus\ZZ/d\ZZ$ with
$d=e(g)\cap[B]$. We call the integer $d\geq0$ the \emph{degree} of $X$.
Bundles of degree $d=0$ are 2-dimensional complex tori. A smooth compact
complex surface $X$ with an elliptic bundle structure of degree $d>0$ is
called a \emph{primary Kodaira surface}. For simplicity, we refer to such
surfaces as Kodaira surfaces.

\medskip Kodaira surfaces have three invariants. The first invariant is
the degree $d>0$. It determines the underlying topological space. The
Universal Coefficient Theorem gives
$$
H^1(X,\ZZ)=\ZZ^3\quadand H^2(X,\ZZ)=\ZZ^4\oplus\ZZ/d\ZZ.
$$
Hence the degree $d$ is also the order of the torsion subgroup of
the N\'eron-Severi group $\NS(X)=\Pic(X)/\Pic^0(X)$.

The second invariant is the $j$-invariant of $B$. It depends only on $X$:
Since $b_1(X)=3$ is odd, $X$ is nonalgebraic. Moreover, since $f:X\ra B$
has connected fibers it is the algebraic reduction and so the fibration
structure does not depend on choices.

The third invariant is an element $\alpha\in\Delta^*$. Here
$\Delta^*= \left\{ z\in\CC\mid 0< |z| <1 \right\}$ is the punctured
unit disk. The Jacobian $\Pic^0_X$ has dimension $h^1(\O_{X})=2$, and
the quotient $\Pic^0_X/\Pic^0_B$ is isomorphic to $\CC^*$. We have
$h^2(\O_{X}(-X_b))=h^0(\O_{X}(X_b))=1$ for each $b\in B$. So the map
on the left in the exact sequence
$$
H^1(X,\O_{X}) \lra H^1(X_b,\O_{X_b}) \lra H^2(X,\O_{X}(-X_b)) \lra
H^2(X,\O_{X}) \lra 0
$$
is surjective. Thus $\Pic^0_X/\Pic^0_B \ra \Pic^0_{X_b}$ is an
epimorphism. By semicontinuity, the kernel equals $\NS(B)=\ZZ$ and is
generated by a well-defined element $\alpha\in\Delta^*$.

\medskip According to \cite{Barth; Peters; Van de Ven 1984}, p.~145, the
principal $E$-bundle $X$ is an associated fibre bundle $X=(L\setminus
0)\times_{\CC^*} \CC^*/\langle\alpha\rangle$ for some line bundle $L\ra B$
of degree $d$. One can check that this $\alpha$ agrees with the invariant
$\alpha$ above. It follows that the isomorphism classes of Kodaira
surfaces correspond bijectively to the triples $(d,j,\alpha)\in
\ZZ_{>0}\times\CC\times\Delta^*$. Moreover, each Kodaira surface with
invariant $(d,j,\alpha)$ is the quotient of a properly discontinuous free
$\ZZ^2$-action on $\CC^*\times\CC^*$ defined by
$$
\Phi(z_1,z_2)= (\beta z_1,z_1^d z_2)\quadand \Psi (z_1,z_2)=
(z_1,\alpha z_2).
$$
Here $\beta\in\Delta^*$ is defined as follows: Consider the
$j$-invariant as a function $j:H\to\CC$ on the upper half plane. It
factors over $\exp:H\to\Delta^*$, $\tau\mapsto \exp(2\pi
\sqrt{-1}\tau)$. Now $\beta=\exp(\tau)$ with $j(\tau)=j$. This
uniformization illustrates Borcea's result on the moduli space of
Kodaira surfaces \cite{Borcea 1984}. Moreover, it suggests the
application of toric geometry for the construction of degenerations.

The following observation will be useful in the sequel:
\begin{lemma}
\label{quotient}
Suppose $X\ra B$ is a principal $E$-bundle of degree $d\geq0$. Let
$G$ be a finite group of order $w$ acting on $X$. If the induced
action on $B$ is free, then $w$ divides $d$, and the quotient $ X/G$
is a principal $E$-bundle of degree $d/w$.
\end{lemma}
\proof Set $X'=X/G$ and $B'=B/G$. Then $B'$ is an elliptic curve, $B\to
B'$ is Galois and the induced fibration $X'\ra B'$ is a principal
$E$-bundle. To calculate its degree $d'$, we use a characterization of
$d$ in terms of $\pi_1(X)$. Since the universal covering of $E$ is
contractible the higher homotopy groups of $E$ vanish. The beginning of
the homotopy sequence of the fibrations $X\to B$, $X'\to B'$ thus leads to
the following diagram of central extensions:
$$
\begin{CD}
  0@>>>\pi_1(E) @>>>\pi_1(X) @>>>\pi_1(B) @>>>0\\
  && @V\id VV @VVV @VVV\\
  0@>>>\pi_1(E) @>>>\pi_1(X') @>>>\pi_1(B') @>>>0\ .
\end{CD}
$$
Let $h_1,h_2\in\pi_1(X')$ map to generators of $\pi_1(B')$. Since
$[\pi_1(B'):\pi_1(B)]=w$ there exists an integral matrix $\binom{a\ b}{c\
d}$ with determinant $w$ such that
$$
g_1=ah_1+bh_2\,,\quad g_2=ch_1+dh_2
$$
map to generators of $\pi_1(B)$. Then $g_1,g_2\in\pi_1(X)\subset
\pi_1(X')$ and $[g_1,g_2]$, $[h_1,h_2]$ are generators for the commutator
subgroups $[\pi_1(X),\pi_1(X)]$ and $[\pi_1(X'),\pi_1(X')]$ respectively.
Now the degrees $d,d'$ being the orders of the torsion subgroups of
$H_1(X,\ZZ) =\pi_1(X)/[\pi_1(X),\pi_1(X)]$ and of $H_1(X',\ZZ)
=\pi_1(X')/[\pi_1(X'),\pi_1(X')]$ they can be expressed as divisibility of
a generator of the commutator subgroup. Write $[h_1,h_2]=d'\cdot x$ for
some primitive $x\in\pi_1(E)$. Then
$$
[g_1,g_2]=[ah_1+bh_2,ch_1+dh_2]=(ad-bc)[h_1,h_2]=wd'\cdot x
$$
shows $d=wd'$ as claimed.\qed

\section{$D$-semistable surfaces with trivial canonical class }
\label{d-semistable}

Our objective is the study of degenerations of smooth Kodaira
surfaces. We consider the following class of singular surfaces:
\begin{definition}
\label{admissible}
A reduced compact complex surface $X$ is called \emph{admissible} if
it is irreducible, has locally normal crossing singularities, is 
$d$-semistable, and
satisfies $K_X=0$.
\end{definition}

The sheaf of first order deformations
$\shT^1_X=\shExt^1(\Omega^1_X,\O_X)$ is supported on $D=\Sing(X)$.
Following Friedman \cite{Friedman 1983}, we call a locally normal
crossing surface $X$ \emph{d-semistable} if $\shT^1_X\simeq \O_D$.
This is a necessary condition for the existence of a global smoothing
with smooth total space.

Smooth admissible surfaces are either K3 surfaces, 2-dimensional
complex tori, or Kodaira surfaces. Throughout this section, $X$ will
be a \emph{singular} admissible surface. We will see that three types
of such surfaces are possible. A finer classification is deferred to
subsequent sections.

\medskip Let $\nu:S\ra X $ be the normalization, $C\subset S$ its
reduced ramification locus, $D\subset X $ the reduced singular locus,
and $\varphi:C\ra D$ the induced morphism. Then $S$ is smooth. The
surface $X$ can be recovered from the commutative diagram
$$
\begin{CD}
  C @>>> S \\
  @V\varphi VV @VV\nu V\\
  D @>>> X,
\end{CD}
$$
which is cartesian and cocartesian. It gives rise to a long exact
Mayer-Vietoris sequence
$$
\ldots \lra H^{p}(X,\O_{X}) \lra H^{p}(S,\O_{S})\oplus
H^{p}(D,\O_{D})\lra H^{p}(C,\O_{C}) \lra \ldots.
$$
The ideals $\O_{S}(-C)\subset\O_{S}$ and $\shI_D= \nu_*(\O_S(-C))
\subset\O_{X} $ are the conductor ideals of the inclusion $\O_{X}
\subset\nu_*(\O_{S}) $. They coincide with the relative dualizing
sheaf
$$
\nu_*(\omega_{S/X}) =\Hom(\nu_*(\O_{S}),\O_{X}).
$$
Hence $K_S=-C$; in particular, the Kodaira dimension is
$\kappa(S)=-\infty$. The Enriques-Kodaira classification of surfaces
(\cite{Barth; Peters; Van de Ven 1984}, Chap.~VI) tells us that $S$ is
either ruled or has $b_1(S)=1$. In the latter case, one says that $S$
is a surface of \emph{class VII}.

Following Deligne and Rapoport \cite{Deligne; Rapoport 1972} we call
the seminormal curve obtained from $\PP^1\times\ZZ/n\ZZ $ by the
relations $(0,i)\sim(\infty,i+1) $, a \textit{N\'eron polygon}.
\begin{lemma}
\label{C}
Each connected component of $C\subset S $ is an elliptic curve or a
N\'eron polygon. For each singular connected component $D'\subset D
$, the number of irreducible components in $\nu^{-1}(D') \subset C$ is
a multiple of 6.
\end{lemma}
\proof The adjunction formula gives $\omega_C=\O_{C} $. Obviously, $C
$ has ordinary nodes. By \cite{Deligne; Rapoport 1972}, Lemma 1.3,
each component of $C$ is an elliptic curve or a N\'eron polygon. Let
$s $ be the number of singularities and $c $ the number of irreducible
components in $D'$. The preimage $\nu^{-1}(D')\subset S $ is a
disjoint union of N\'eron polygons. It has $2c $ irreducible
components and $3s $ singularities. For N\'eron polygons, these
numbers coincide. \qed

We will use below the following {\em triple point formula} for
$d$-semistable normal crossing surfaces:
\begin{lemma}
\label{triple}
(\cite{Persson 1977}, Cor.2.4.2) Let $D'\subset D$ be an irreducible
component and $C'_1\cup C'_2\subset C$ its preimage in $S$. Then
$-(C'_1)^2-(C'_2)^2$ is the number of triple points of $D'$.
\end{lemma}

Let $g:S\ra S'$ be a minimal model and $C'=g(C)$. Then $K_{S'}=-C'$,
so $C'\subset S'$ has ordinary nodes, and $g:S\ra S'$ is a sequence of
blowing-ups with centers over $\Sing(C')$. Recall that $S$ is a
\emph{Hopf surface} if its universal covering space is isomorphic to
$\CC^2\setminus 0$.
\begin{proposition}
\label{Hopf}
Suppose $S$ is nonalgebraic. Then $S$ is a Hopf surface. The curve
$C=C_1\cup C_2 $ is a disjoint union of elliptic curves, and $X $ is
obtained by gluing them together.
\end{proposition}
\proof Since $S$ has algebraic dimension $a(S)<2$, no curve on $S$ has
positive selfintersection. Suppose $C$ contains a N\'eron polygon
$C'=C_1\cup\ldots\cup C_m$. As $X$ is $d$-semistable, the triple point
formula Lemma~\ref{triple} implies $C_i^2=-1$. So $(C')^2=\sum C_i^2
+\sum_{i\neq j} C_i\cdot C_j = -m + 2m>0$, contradiction. Hence $C$
is the disjoint union of elliptic curves. Consider the exact sequence
$$
H^0(S,\O_S) \lra H^0(C,\O_C) \lra H^1(S,\omega_S).
$$
By the Kodaira classification, $b_1(S)=1$, consequently
$h^1(\O_S)=1$, thus $C$ has at most two components. Suppose
$\varphi:C\ra D$ induces a bijection of irreducible components. Then
each component of $C$ double-covers its image in $D$. So the map on
the left in the exact sequence
$$
H^1(S,\O_{S})\oplus H^1(D,\O_{D}) \lra H^1(C,\O_{C}) \lra
H^2(X,\O_{X})\lra 0
$$
is surjective, contradicting $K_X=0$. Consequently $C=C_1\cup C_2$
consists of two elliptic curves which are glued together in $X$. Let
$n\geq0$ be the number of blowing-ups in $S\ra S'$. The normal bundle
of $C$ has degree $K_S^2=K_{S'}^2-n=-n$. Again by the triple point
formula it follows $n=0$. Hence $S$ is minimal. By \cite{Nishiguchi
  1988}, Prop.~3.1, $S$ is a Hopf or Inoue surface. The latter case is
impossible here because Inoue surfaces contain at most one elliptic
curve (\cite{Nishiguchi 1988}, Prop.~1.5).\qed
\begin{proposition}
\label{nonrational algebraic}
Suppose $S$ is nonrational algebraic. Then $S$ is a $\PP^1$-bundle
over an elliptic curve $B$. The curve $C=C_1\cup C_2$ is the disjoint
union of two sections, and $X$ is obtained by gluing them together.
\end{proposition}
\proof By the Enriques classification, the minimal model $S'$ is a
$\PP^1$-bundle over a nonrational curve $B$. By L\"uroth's Theorem,
$C$ does not contain N\'eron polygons, hence is the disjoint union of
elliptic curves. The Hurwitz formula implies that $B $ is also an
elliptic curve. Since $K_S$ has degree $-2$ on the ruling, we infer
$C=-K_S$ has at most two irreducible components. As in the preceding
proof, we conclude that $C$ consists of two components, which are
identified in $X$. Let $n\geq0$ be the number of blowing-ups in $S\ra
S'$. The normal bundle of $C$ has degree $K_S^2= K_{S'}^2-n =-n$.
Since $X$ is $d$-semistable, $n=0$ follows. Hence $S$ is a
$\PP^1$-bundle. \qed

\begin{proposition}
\label{rational}
Suppose $S$ is a rational surface. Then it is the blowing-up of a
Hirzebruch surface in two points $P_1,P_2$ in disjoint fibers
$F_1,F_2$. The curve $C$ is a N\'eron 6-gon, consisting of the two
exceptional divisors, the strict transforms of $F_1,F_2$, and the
strict transforms of two disjoint sections whose union contains $P_1$
and $P_2$. The surface $X$ is obtained by identifying pairs of
irreducible components in $C$.
\end{proposition}
\proof The minimal model $S'$ is either $\PP^2$ or a Hirzebruch
surface of degree $e\neq 1$. Since each representative of $-K_{S'}$
is connected the curve $C$ must be connected. Suppose $C$ is an
elliptic curve. Then the map on the left in the exact sequence
$$
H^1(S,\O_{S})\oplus H^1(D,\O_{D}) \lra H^1(C,\O_{C}) \lra
H^2(X,\O_{X})\lra 0
$$
is surjective, contradicting $K_X=0$. Consequently,
$C=C_1\cup\ldots \cup C_m$ is a N\'eron polygon. The triple point
formula gives $\sum C_i^2=-m$, hence
$$
K_S^2=\sum C_i^2 +\sum_{i\neq j} C_i\cdot C_j = -m+2m =m.
$$
Let $n\geq 0$ be the number of blowing-ups in $S\ra S'$, so $K_S^2
= K_{S'}^2-n$. For $S'=\PP^2$ this gives $m=9-n$. If $S'$ is a
Hirzebruch surface, $m=8-n$ holds instead. According to Lemma
\ref{C}, the natural number $m$ is a multiple of $6$. The only
possibilities left are $S'=\PP^2$ and $C'$ a N\'eron 3-gon, or $S'$ a
Hirzebruch surface and $C'$ a N\'eron 4-gon. From this it easily
follows that $X$ is obtained from a blowing-up of a Hirzebruch surface
as stated. \qed
\begin{remark}
\label{numerically admissible}
For the results of this section one can weaken the hypothesis that $X$
is $d$-semistable. It suffices to assume that the invertible
$\O_D$-module $\shT^1_X$ is numerically trivial. One might call such
surfaces \emph{numerically d-semistable}. They will occur as locally
trivial deformations of $d$-semistable surfaces.
\end{remark}

\section{Hopf surfaces}

In this section we analyze the geometry and the deformations of
admissible surfaces $X$ whose normalization $S$ is
\emph{nonalgebraic}. By Proposition~\ref{Hopf}, such $S$ is a Hopf
surface containing a union of two disjoint isomorphic elliptic curves
$C=C_1\cup C_2$.

As a topological space, Hopf surfaces are certain fibre bundles over
the 1-sphere, whose fibres are quotients of the 3-sphere (see
\cite{Kato 1975}, Thm.~9, and \cite{Kato 1989}). By Hartog's
Theorem, the action of $\pi_1(S) $ on the universal covering space
$\CC^2\setminus 0 $ extends to $\CC^2 $, fixing the origin. We call a
Hopf surface \emph{diagonizable} if $\pi_1(S)$ is contained in the
maximal torus $\CC^*\times\CC^*\subset\GL_2(\CC)$ up to conjugacy
inside the group of all biholomorphic automorphisms of $\CC^2 $ fixing
the origin. This is precisely the class of Hopf surfaces we are
interested in:
\begin{proposition}
\label{diagonizable}
A Hopf surface $S $ is diagonizable if and only if it contains two
disjoint elliptic curves $C_1,C_2\subset S$ with $-K_S=C_1+C_2 $.
Moreover, a diagonizable Hopf surface has fundamental group
$\pi_1(S)=\ZZ\oplus\ZZ/ n\ZZ$.
\end{proposition}
\proof Suppose $\pi_1(S)\subset \GL_2(\CC) $ consists of diagonal
matrices. Let $C_1,C_2\subset S $ be the images of $\CC^*\times 0$
and $ 0\times\CC^* $, respectively. The invariant meromorphic 2-form
$dz_1/z_1\wedge dz_2/z_2 $ on $\CC^2 $ yields $-K_S=C_1+C_2 $.

Conversely, assume the condition $-K_S=C_1+C_2$. An element
$\alpha\in\pi_1(S)$ is called a \emph{contraction} if
$\lim_{n\ra\infty}\alpha^n(U)= \left\{ 0 \right\}$ holds for the unit
ball $U$ in $\CC^2$. According to a classical result \cite{Sternberg
  1957}, equation 44, in suitable coordinates a contraction takes the
form
$$
\alpha(z_1,z_2) = (\alpha_1 z_1+\lambda z_2^m,\alpha_2 z_2)
$$
with $ 0 < |\alpha_1| \leq |\alpha_2|<1$ and
$\lambda(\alpha_1-\alpha_2^m)=0$. We claim that each contraction
$\alpha\in\pi_1(S)$ has $\lambda=0$. Suppose not. The quotient of
$\CC^2\setminus 0$ by $\langle\alpha\rangle$ is a primary Hopf surface
and defines a finite \'etale covering $f:S'\ra S$. Let $C'\subset S'$
be the image of $\CC^*\times 0$. Then \cite{Kodaira 1966}, p.~696,
gives $-K_{S'}=(m+1)C'$. On the other hand,
$$
-K_{S'}=-f^*(K_S)=f^{-1}(C_1)+f^{-1}(C_2)
$$
holds. Consequently, a nonempty subcurve of $(m+1)C'$ is base
point free, so $S'$ is elliptic. Now according to \cite{Kodaira
  1966}, Theorem 31, ellipticity of $S'=(\CC^2\setminus0)/(\alpha)$
implies $\lambda=0$, contradiction.

Next we claim that $\pi_1(S)$ is contained in $\GL_2(\CC)$, at least
up to conjugacy. Suppose not. A nonlinear fundamental group
$\pi_1(S)$ is necessarily abelian (\cite{Kodaira 1966}, Thm.~47).
By \cite{Kato 1975}, p.~231, there is a contraction
$\alpha\in\pi_1(S)$ with $\lambda\neq 0$, contradiction.

Now we can assume $\pi_1(S)\subset\GL_2(\CC)$. Write $\pi_1(S)=
G\cdot Z$ as a semidirect product, where $G=\left\{ \gamma\in\pi_1(S)
  \mid |\det(\gamma)| = 1 \right\}$, and $Z\simeq\ZZ$ is generated by
a contraction $\alpha=(\alpha_1,\alpha_2)$, acting diagonally. If
$\alpha_1\neq\alpha_2$ then \cite{Wehler 1981}, p.~24, ensures that
$S$ is diagonizable. It remains to treat the case
$\alpha_1=\alpha_2$. Then $Z$ is central and we only have to show that
$G$ is abelian. Let $S'$ be the quotient of $\CC^2\setminus 0$ by
$\langle\alpha\rangle$. The canonical projection $\CC^2\setminus 0\ra
\PP^1$ induces an elliptic bundle $S'\ra\PP^1$. This defines an
elliptic structure $g:S=S'/G\lra \PP^1/G$. Suppose $\PP^1\to\PP^1/G$
has more than $2$ ramification points. Let $F_b$ be the reduced fiber
over $g\in\PP^1/G$ with multiplicity $m_b$ and $C_i=F_{b_i}$; the
canonical bundle formula yields
\begin{eqnarray*}
0&=&C_1+C_2+ K_S
\ =\ F_{b_1}+F_{b_2}-2F+\sum_{b\in\PP^1/G}(m_b-1)F_b\\
&=&(1-m_{b_1})F_{b_1}+(1-m_{b_2})F_{b_2}
+\sum_{b\in\PP^1/G}(m_b-1)F_b\ >\ 0\,,
\end{eqnarray*}
contradiction. Consequently there are only 2 ramification points.
Therefore $G$ is a subgroup of $\CC^*$, hence abelian. This finishes
the proof of the equivalence.

As for the fundamental group, it must have rank 1 for $b_1=1$. If
there were more than one generator needed for the torsion part the
action could not be free on $(\CC^*\times 0) \cup ( 0\times\CC^*)$.
\qed

\medskip Suppose that $S$ is a diagonizable Hopf surface with
$\pi_1(S)\simeq \ZZ\oplus\ZZ/n\ZZ \subset \GL_2(\CC) $ consisting of
diagonal matrices. The free part of $\pi_1(S)$ is generated by a
contraction $\alpha=(\alpha_1,\alpha_2)$ with $\alpha_1,
\alpha_2\in\Delta^*$. The torsion part is generated by a pair
$\zeta=(\zeta_1,\zeta_2)$ of primitive $n$-th roots of unity. Note
that they must be primitive, again since the action is free on the
coordinate axes minus the origin. Choose $\tau_i\in \CC$ and
$n_i\in\ZZ$ with $\alpha_i =\exp(2\pi\sqrt{-1}\tau_i)$ and
$\zeta_i=\exp(2\pi\sqrt{-1}n_i/n)$. Consider the lattice
$$
\Lambda_i= \left< \tau_i, n_i/n, 1\right> = \ZZ\cdot\tau_i + \ZZ\cdot 1/n
\subset \CC.
$$
The images $C_1,C_2\subset S $ of $\CC^*\times 0$ and $0\times\CC^*$
take the form
$$
C_i= \CC^*/\left<\alpha_i,\zeta_i\right> = \CC / \Lambda_i.
$$
Without loss of generality we can assume $C=C_1\cup C_2$. In fact,
if $S$ is not elliptic then $C_1,C_2$ are the only curves on $S$ at
all. Otherwise the elliptic fibration comes from a rational function
of the form $z_1^a/z_2^b$ on $\CC^2\setminus 0$. Then $|ab|\neq 1$ iff
the elliptic fibration $S\to\PP^1$ has multiple fibers. In this case
the canonical bundle formula implies $C=C_1\cup C_2$. Otherwise we can
apply a linear coordinate change on $\CC^2$ to achieve that the
preimages of the components of $C$ are $\CC^*\times 0$ and
$0\times\CC^*$.

Now assume there is an isomorphism $\psi:C_2\ra C_1$. Set $D=C_1$,
$C=C_1\cup C_2$, and $\varphi=\id\cup\psi:C\ra D$. The cocartesian
diagram
$$
\begin{CD}
  C @>>> S \\
  @V \varphi VV @VV \nu V \\ D @>>> X
\end{CD}
$$
defines a normal crossing surface $X$. Since each translation of
the elliptic curve $C_2$ extends to an automorphism of $S$ fixing
$C_1$, we can assume that $\psi:C_2\ra C_1$ respects the origin.
Hence $\psi$ is defined by a homothety $\mu:\CC\ra\CC$ with
$\mu\Lambda_2=\Lambda_1$. In other words, there is a matrix
\begin{equation}
\label{homothety}
\begin{aligned}
\begin{pmatrix}
  a & c \\ b & d
\end{pmatrix}
\in\SL_2(\ZZ)
\end{aligned}
\qquad\text{with}\qquad
\begin{aligned}
  \mu\tau_2 = a\tau_1 + b/n \\
  \mu/n = c\tau_1 + d/n.
\end{aligned}
\end{equation}
The following two propositions describe the properties of the surface
$X$ in terms of $\pi_1(S)\subset \GL_2(\CC)$ and $\mu\in \CC$.
\begin{proposition}
\label{untwisted Hopf}
The condition $K_X=0$ holds if and only if $a=d=1$ and $c=0$.
\end{proposition}
\proof As $\nu^*(K_X)=0$ in any case $K_X$ is numerically trivial.
The condition $K_X=0$ is thus equivalent to $H^2(X,\O_{X})\neq 0$.
Assume $K_X=0$. Then the map on the left in the exact sequence
$$
H^1(S,\O_{S})\oplus H^1(D,\O_{D}) \lra H^1(C,\O_{C}) \lra
H^2(X,\O_{X}) \lra 0
$$
is not surjective. Hence $H^1(\O_{D})$ and $H^1(\O_{S})$ have the
same image in $H^1(\O_{C})$. According to the commutative diagram
$$
\begin{CD}
  H^1(S,\O_{S}) @>>> H^1(C,\O_{C}) @< \varphi^*<< H^1(D,\O_{D}) \\
  @V\text{bij} VV @VV\text{inj}V @VV\text{inj}V\\
  H^1(S,\CC) @>>> H^1(C,\CC) @< \varphi^*<< H^1(D,\CC)
\end{CD}
$$
from Hodge theory, the image of $H^1(D,\CC)$ in $H^1(C,\CC)$
contains the image of $H^1(S,\CC)$. So the composition
$$
H^1(S,\CC) \lra H^1(C,\CC) \stackrel{(\id,\psi^*)^{-1}}{\lra}
H^1(D,\CC)\oplus H^1(D,\CC)
$$
factors over the diagonal $H^1(D,\CC)\subset H^1(D,\CC)^{\oplus
  2}$. Dually, the composition
$$
H_1(D,\CC)\oplus H_1(D,\CC) \stackrel{(\id,\psi_*)^{-1}}{\lra}
H_1(C,\CC)\lra H_1(S,\CC)
$$
factors over the addition map $H_1(D,\CC)^{\oplus 2}\ra
H_1(D,\CC)$. In other words, the composition
$$
\Lambda_1 \stackrel{(\id, -\id)}{\lra} \Lambda_1\oplus\Lambda_1
\stackrel{\id\oplus\binom{\mbox{}\ d\ -c}{-b\ \ a} }{\lra}
\Lambda_1\oplus\Lambda_2 \stackrel{(1,0,1,0)}{\lra} \ZZ
$$
is zero. It is described by $(1-d,c) $; thus $d=1$, $c=0 $, and in
turn $a=1 $. Hence the condition is necessary. The converse is shown
in a similar way. \qed
\begin{remark}
  Inserting the above conditions into equation \ref{homothety}, one
  sees that the homothety $\mu:\CC\ra\CC$ must be the identity, and
  then $\tau_2=\tau_1 +b/n$. In particular $S$ is elliptically
  fibered.
\end{remark}
\begin{proposition}
\label{d-semistable Hopf}
Suppose $K_X=0$. Then $X$ is $d$-semistable if and only if the
congruences $(n_1-n_2)^2 \equiv 0 $ and $b(n_1-n_2)\equiv 0$ modulo
$n$ hold.
\end{proposition}
\proof Let $ N_i$ be the normal bundle of $C_i\subset S $. We can
identify $\shT^1_X $ with $N_2\otimes \psi^*(N_1) $. The bundle $N_1$
can be obtained as quotient of the normal bundle of $\CC^*\times 0$ in
$\CC^*\times\CC^*$, which is trivial, by the induced action of
$\pi_1(S)=\ZZ\oplus\ZZ/n \ZZ$: $N_1=\CC^*\times\CC^*/ (\alpha,\zeta)$
with
$$
\alpha\cdot (z_1,v) = (\alpha_1z_1, \alpha_2v) \quadand
\zeta\cdot(z_1,v) = (\zeta_1z_1, \zeta_2v).
$$
After pull back along $\CC\ra \CC^* ,z\mapsto\exp(2\pi\sqrt{-1}z)$,
the bundle $N_1$ can also be seen as the quotient of the trivial
bundle $\CC\times\CC $ by $\pi_1(C_1)=\Lambda_1 $ acting via
$$
\tau_1\cdot(w,v) = (w+\tau_1,\alpha_2v),\quad n_1/n\cdot(w,v) =
(w+n_1/n, \zeta_2 v),\quad 1\cdot (w,v) =(w+1,v).
$$
Choose $m_1$ with $m_1n_1\equiv1$ modulo $n$. Then
$1/n\in\Lambda_1 $ acts via $1/n\cdot(w,v)=(w+1/n,\zeta_2^{m_1}v) $.
An analogous situation holds for $N_2 $. Inserting $\Psi=\binom{1\
  0}{b\ 1}$ we obtain that $\shT^1_X=N_2\otimes\psi^*(N_1) $ is the
quotient of $\CC\times\CC $ by the action
$$
\tau_2\cdot(w,v)=(\tau_2+w,\alpha_1\alpha_2\zeta_2^{bm_1}v)
\quadand 1/n\cdot(w,v)= (1/n + w, \zeta_1^{m_2}\zeta_2^{m_1} v).
$$
The isomorphism class of $\shT^1_X=N_2\otimes\psi^*(N_1) $ depends
only on the element
$$
\tau_2\longmapsto\alpha_1\alpha_2\zeta_2^{bm_1},\quad 1/n
\longmapsto \zeta_1^{m_2}\zeta_2^{m_1} $$
in
$\Hom_\ZZ(\Lambda_2,\CC^*) $. To check for triviality of $\shT^1_X $,
consider the commutative diagram
$$
\begin{CD}
  \Hom_\ZZ(\Lambda_2,\CC) @>\pr>> \Hom_\CC(\CC,\overline{\CC})
  @>\simeq>> H^1(C_2,\O_{C_2})\\ @V\exp VV @. @VV\exp V\\
  \Hom_\ZZ(\Lambda_2,\CC^*) @>>> \Pic(C_2) @<<\simeq<
  H^1(C_2,\O_{C_2}^*),
\end{CD}
$$
explained in \cite{Mumford 1970}, Sect.\ I.2. Here $\pr $ is the
projection of $\Hom_\ZZ(\Lambda_2,\CC)\simeq \Hom_\RR(\CC,\CC) $ onto
the $\CC $-antilinear homomorphisms $\Hom_\CC(\CC,\overline{\CC}) $.
Lifting the homomorphism defining $N_2\otimes\varphi^*(N_1) $ from
$\Hom_\ZZ(\Lambda_2,\CC^*) $ to $\Hom_\ZZ(\Lambda_2,\CC) $, we obtain
$$
\tau_2\longmapsto \tau_1+\tau_2+ \frac{bm_1n_2}{n},\quad
1/n\longmapsto \frac{m_2n_1 + m_1n_2}{n}.
$$
Using the diagram, we infer that $X $ is $d$-semistable if and only
if this homomorphism is $\CC $-linear up to an integral homomorphism.
The latter condition is equivalent to the equation
$$
(\tau_1+\tau_2+\frac{bm_1n_2}{n} + e) - \tau_2n(\frac{m_2n_1 +
  m_1n_2}{n} +f)=0
$$
for certain integers $e,f $. Now proceed as follows: Substitute
$\tau_1=\tau_2-b/n$; the coefficients of $\tau_1$ and $1$ give two
equations; since $e$ and $f$ are variable this leads to two
congruences modulo $n$; finally use $n_1m_1\equiv 1$, $n_2m_2\equiv 1$
modulo $n$ to deduce the stated congruences. \qed

\medskip We seek a coordinate free description of $X$. We call an
automorphism of a N\'eron 1-gon a \emph{rotation} if the induced
action on the normalization $\PP^1$ fixes the branch points (compare
\cite{Deligne; Rapoport 1972}, Sect.~3.6).
\begin{proposition}
\label{cover Hopf}
Let $X$ be a complex surface. The following are equivalent:
\renewcommand{\labelenumi}{(\roman{enumi})}
\begin{enumerate}
\item $X$ is admissible with nonalgebraic normalization $S$.
\item There is an elliptic principal bundle $X'\ra B'$ of degree $e>0$
  over the N\'eron 1-gon $B'$ and a (cyclic) Galois covering $g:X'\ra
  X$ such that the Galois group $G$ acts effectively on $B'$ via
  rotations.
\end{enumerate}
\end{proposition}
\proof Suppose (ii) holds. The normalization $S'$ of $X'$ is an
elliptic principal bundle of degree $e>0$ over $\PP^1$, hence
nonalgebraic. Since $g:X'\ra X$ induces a finite morphism $g':S'\ra
S$, the surface $S$ is nonalgebraic as well. It is easy to check
$K_X=0$ using the assumption that $G$ acts via rotations.

Conversely, assume (i). As we have seen, $S$ is a diagonizable Hopf
surface. The exact sequence
$$
0 \lra \pi_1(S') \lra \pi_1(S) \lra \PGL_2(\CC)
$$
determines another diagonizable Hopf surface $S'$ and a Galois
covering $S'\ra S$. The Galois group $G=\pi_1(S)/\pi_1(S')$ is cyclic
of order $m= n/\gcd(n,n_1-n_2,b)$, and the order of the torsion
subgroup in $\pi_1(S')$ is $n'=\gcd(n,n_1-n_2)$, as simple
computations show.

The projection $\CC^2\setminus 0\ra \PP^1$ induces an elliptic
principal bundle $f':S'\ra \PP^1$. The Galois group $G$ acts
effectively on $\PP^1$ and fixes $0,\infty$. Let $C_1',C_2'\subset
S'$ be the fibres over the fixed points. The quotient $\PP^1/G$ is a
projective line, and we obtain an induced elliptic structure $f:S\ra
\PP^1/G$. Note that $f$ has two multiple fibres of order $m$, whose
reductions are $C_1,C_2$. Next, the canonical isomorphism
$\psi':C'_2\ra C_1'$ gives a normal crossing surface $X'$. Clearly,
it is an elliptic principal bundle of degree $e=n'>0$ over the N\'eron
1-gon $B'$ obtained from $\PP^1$ by the relation $0\sim\infty$. The
congruences $(n_1-n_2)^2\equiv0\equiv b(n_1-n_2)$ modulo $n$ are
equivalent to $n'(n_1-n_2)\equiv0\equiv n'b$. A straightforward
calculation shows that this ensures that the free $G$-action on $S'$
descends to a (free) action on $X'$. The quotient is $X=X'/G$. \qed
\begin{remark}
\label{classify Hopf}
At this point we are in position to classify the admissible surfaces
with nonalgebraic normalization. The elliptic principal bundles $X'\to
B'$ are in one-to-one correspondence with pairs $(\alpha,e)\in
\Delta^*\times \NN$, while up to isomorphism the action of $G$ on $B'$
only depends on $|G|$. A further discrete invariant belongs to
possibly non-isomorphic lifts of the $G$-action to $X'$.
\end{remark}
We call the invariants $e$ and $w=|G|$ of $X$ the
\emph{degree} and the \emph{warp} of $X$ respectively. The congruences
in Proposition~\ref{d-semistable Hopf} imply that $w$ divides $e$.
This is important for the construction of smoothings:
\begin{theorem}
\label{smoothing Hopf}
Suppose $X$ is an admissible surface (Definition~\ref{admissible})
with nonalgebraic normalization, of degree $e$ and warp $w$. Then $X$
deforms to a smooth Kodaira surface of degree $e/w$.
\end{theorem}
\proof Let $X'\ra X$ be the Galois covering from Proposition
\ref{cover Hopf}. It suffices to construct a $G$-equivariant
deformation of $X'$.

The deformation of $X'$ will be obtained as infinite quotient of a
toric variety belonging to an infinite fan. Our construction fibers
over a similar construction of a smoothing of the N\'eron 1-gon due to
Mumford (\cite{Ash; Mumford; Rapoport; Tai 1975}, Ch.I,4) that we now
recall for the reader's convenience. Mumford considered the fan
generated by the cones
$$
\bar\sigma_m=\Big\langle(m,1),(m+1,1)\Big\rangle,\quad m\in\ZZ.
$$
There is a $\ZZ$-action on the associated toric variety $W$
belonging to the linear \mbox{(right-)} action by $\binom{1\ 0}{1\ 1}$
on the fan. Let $q:W\to\CC$ be the morphism coming from the projection
$\pr_2:\ZZ^2\to\ZZ$. Note that the central fiber is an infinite chain
of rational curves (a N\'eron $\infty$-gon), while the general fiber
is isomorphic to $\CC^*$. Now $\ZZ$ acts fiberwise and the action is
proper and free over the preimage of $\Delta\subset\CC$. Moreover,
$\ZZ$ acts transitively on the set of irreducible components of the
central fiber. The quotient is the desired deformation of the N\'eron
1-gon.

For our purpose, consider the infinite fan in $\ZZ^3$ generated by the
two-dimensional cones
$$
\sigma_m = \Big\langle(m,e\cdot{\textstyle \binom{m}{2}}, 1),
(m+1,e\cdot{\textstyle \binom{m+1}{2}}, 1)\Big\rangle,\quad m\in\ZZ.
$$
Let $V$ be the corresponding 3-dimensional smooth toric variety and
$V\to \CC$ the toric morphism belonging to the projection
$\pr_3:\ZZ^3\ra\ZZ$. Now the special fibre $V_0$ is isomorphic to a
$\CC^*$-bundle over the N\'eron $\infty$-gon. As neighbouring cones
$\sigma_m$, $\sigma_{m+1}$ are not coplanar, its degree over each
irreducible component of the N\'eron $\infty$-gon is nonzero. A
simple computation shows that its degree is $e$. This time the
$\ZZ$-action on the fan is given by the automorphism
$$
\Phi=
\begin{pmatrix}
  1&e&0 \\
  0& 1&0\\
  1&0&1
\end{pmatrix}
\in\SL_3(\ZZ).
$$
Again the induced $\ZZ$-action on the preimage $U\subset V$ of
$\Delta\subset \CC$ is proper and free. Choose $\alpha\in\Delta^*$
such that the fibres of $X'\ra B'$ are isomorphic to
$\CC^*/\langle\alpha^w\rangle$. Let $\Psi:V\ra V$ be the automorphism
extending the action of $(1,\alpha,1)\in(\CC^*)^3$. Then
$\foX'=U/\langle \Phi,\Psi^w\rangle$ is a smooth 3-fold, endowed with
a projection $\foX'\ra\CC$. The general fibres $\foX_t$,
$t\in\Delta^*$, are smooth Kodaira surfaces of degree $e$ with fibres
$\CC^*/\langle\alpha^w\rangle$ and basis $\CC^*/\langle t\rangle$.
The special fibre $\foX'_0$ is isomorphic to $X'$.

Finally we extend the $G$-action on $X'$ to $\foX'$. The automorphism
$\Psi:V\ra V$ descends to a free $G$-action on $\foX'$. Replacing
$\alpha$ by another primitive $w$-th root of $\alpha^w$ if necessary,
the induced action on $\foX'_0$ coincides with the given action on
$X'$. Consequently, $\foX=\foX'/G$ is the desired smoothing of
$X=\foX_0$. The action of $\Psi$ on $B'$ is free. Therefore, according
to Lemma \ref{quotient}, the general fibres $\foX_t$ are smooth
Kodaira surfaces of degree $e/w$. \qed

\medskip The next task is to determine the versal deformation of $X$.
We have to calculate the relevant cohomology groups:
\begin{lemma}
\label{cohomology Hopf}
For a locally normal crossing surface $X$ with $K_X=0$ and
non-algebraic normalization it holds $h^0(\Theta_X)=1 $,
$h^1(\Theta_X)=1 $ and $h^2(\Theta_X)=0 $.
\end{lemma}
\proof The Ext spectral sequence together with a local computation
gives
$$
H^p(X,\Theta_X)=\Ext^p(\Omega^1_X/\tau_X,\O_{X}),
$$
see \cite{Friedman 1983}, page 88ff. Here $\tau_X\subset \Omega_X^1 $
denotes the torsion subsheaf. In view of Serre duality
$$
\Ext^p(\Omega^1_X/\tau_X,\O_{X})^\vee \simeq
H^{2-p}(X,\Omega^1_X/\tau_X)
$$
it thus suffices to compute the cohomology of $\Omega^1_X/\tau_X$.
Consider the exact sequence
$$
0 \lra \Omega^1_X/\tau_X \lra \Omega_S^1 \lra \Omega^1_D \lra 0
$$
of $\O_X$-modules (\cite{Friedman 1983}, Prop.~1.5). Note
that the map on the right is an \emph{alternating sum}. The inclusion
$H^0(X,\Omega^1_X/\tau_X) \subset H^0(S,\Omega_S^1) =0$ yields $ h^2(
\Theta_X)=0$. Moreover,
$$
0 \lra H^0(D,\Omega^1_D) \lra H^1(X,\Omega^1_X/\tau_X) \lra
H^1(S,\Omega^1_S)
$$
and $h^{1,1}(S)=0$ gives $h^1(\Theta_X)=1 $. With $b_3(S)=1 $,
$h^1(\omega_S)=h^1(\O_{S})=1 $ and degeneracy of the Fr\"ohlicher spectral
sequence for smooth compact surfaces (\cite{Barth; Peters; Van de Ven
1984}, IV, Thm.~2.7) we have $h^{1,2}(S)=0 $. Now the exact sequence
$$
0 \lra H^1(D,\Omega^1_D) \lra H^2(X,\Omega^1_X/\tau_X) \lra
H^2(S,\Omega^1_S)
$$
gives $h^0(\Theta_X)=1 $. \qed
\begin{proposition}
\label{T^i Hopf}
Let $X$ be as in Theorem~\ref{smoothing Hopf}. Then $\dim T^0_X=1$,
$\dim T^1_X=2$ and $\dim T^2_X=1$.
\end{proposition}
\proof Since $X$ is locally a complete intersection the
$E_2$ term of the spectral sequence $E_2^{p,q}=
H^p(\shT^q_X) \Rightarrow T^{p+q}_X$ has at most one non-trivial
differential, which is $H^0(\shT^1_X)\to H^2(\Theta_X)$. The
previous proposition shows first degeneracy at $E_2$ level
and in turn gives the stated values for $\dim T^i_X$.
\qed
\begin{theorem}
\label{versal Hopf}
  $X $ is an admissible surface with nonalgebraic normalization,
of degree $e$ and warp $w$. Then the semiuniversal deformation
$p:\foX\ra V $ of $X $ has a smooth, 2-dimensional base $V$. The
locally trivial deformations are parameterized by a smooth curve
$V'\subset V $, and $V\setminus V' $ corresponds to smooth Kodaira
surfaces of degree $e/w$.
\end{theorem}
\proof Since $h^1(\Theta_X)=1 $ and $h^2(\Theta_X)=0 $, the locally
trivial deformations are unobstructed, and $V' $ is a smooth curve.
Since $\dim T^1_X=2 $ and since $X $ deforms to smooth Kodaira
surfaces of degree $e/w$, which move in a 2-dimensional family, the
base $ V$ is a smooth surface. We saw in the proof of the previous
proposition that the restriction map $T^1\ra H^0(X,\shT^1_X)$ is
surjective. According to \cite{Friedman 1983}, Proposition 2.5, the
total space $\foX$ is smooth. Now Sard's Lemma implies that the
projection $\foX\ra V$ is smooth over $V\setminus V'$, at least after
shrinking $V$.
\qed

\begin{remark}
The referee pointed out that $V'$ should have an interpretation as
versal deformation space of the elliptic curve $D\subset X$. This is
indeed the case: Since the restriction of the Kodaira-Spencer map
to $T_{V'}$ generates $H^1(\Theta_X)$ it suffices to show that the
composition
\begin{eqnarray} \label{KS of Hopf}
H^1(\Theta_X)\lra H^1(\Theta_X\otimes\O_D)\lra H^1(\Theta_D)
\end{eqnarray}
is surjective. This statement is stable under \'etale covers, so we
may assume that there is an elliptic fibration $p:X\to B$ over the
N\'eron 1-gon (Proposition~\ref{cover Hopf}). The double curve $D$ is
the fiber over the node of $B$, and by relative duality
$R^1p_*(\Theta_{X/B})= \O_B$. The base change theorem implies
that the restriction map $R^1p_*(\Theta_{X/B})\to H^1(\Theta_D)$ is
surjective. On global sections this map is nothing but the
composition of 
$$
\CC=H^0(B,R^1p_*(\Theta_{X/B}))\lra H^1(X,\Theta_X)
$$
with (\ref{KS of Hopf}). Therefore the latter map is surjective too.
\end{remark}

\section{Ruled surfaces over elliptic curves}

In this section we study admissible surfaces $X $ whose normalization
$S $ is \emph{algebraic and nonrational}, as in Proposition
\ref{nonrational algebraic}.

Let $B'$ be an elliptic curve and $f:S\ra B'$ a $\PP^1$-bundle with
two disjoint sections $ C_1, C_2\subset S$. Put $C=C_1\cup C_2$ and
$D=C_1$. Let $\psi:C_2\ra C_1$ be an isomorphism and
$\varphi=\id\cup\psi:C\ra D$ the induced double covering. The
cocartesian diagram
$$
\begin{CD}
  C @>>> S \\
  @V\varphi VV @VV\nu V\\
  D @>>> X
\end{CD}
$$
defines a normal crossing surface $X$. By abuse of notation we
write $\psi$ also for the induced automorphism $(f|_{C_1})\circ
\psi\circ (f|_{C_2})^{-1}$ of $B'$.
\begin{proposition}
\label{untwisted elliptic}
We have $K_X=0$ if and only if $\psi:B'\ra B'$ is a translation.
\end{proposition}
\proof As in the proof of Proposition~\ref{untwisted Hopf}, $K_X$ is
trivial iff the images of $H^1(S,\O_S)$ and $H^1(D,\O_D)$ in
$H^1(C,\O_C)$ coincide. Notice that the map $\phi^*:H^1(\O_{D})\ra
H^1(\O_{C}) \simeq H^1(\O_D)\oplus H^1(\O_D) $ is the diagonal
embedding. It does not depend on $\phi$.

Suppose $\psi$ is a translation. Then the images of $H^1(\O_{S}) $
and $H^1(\O_{D})$ both agree with $(f|_C)^*H^1(\O_{B'})$.

Conversely, assume that $\psi$ is not a translation. Then it acts on
$H^1(\O_{E})$ as a nontrivial root of unity. Consequently,
$H^1(\O_{S}) $ and $H^1(\O_{C}) $ have different images in
$H^1(\O_{D}) $. \qed

\medskip The next task is to calculate the sheaf of first order
deformations $\shT^1_X$. We call the degree $e=-\min \left\{ A^2 \mid
  A\subset S \;\text{ a section} \right\}$ of $S$ also the
\emph{degree} of $X$. In our case $e\geq 0$, since $S\ra B'$ has two
disjoint sections. Let $w$ be the order, possibly $0$, of the
automorphism $\psi:B'\ra B'$. Call it the \emph{warp} of $X$.

\begin{proposition}
\label{d-semistable elliptic}
Suppose $K_X=0$. Then $X$ is $d$-semistable if and only if the warp $w$
divides the degree $e$.
\end{proposition}
\proof We identify $C_1,C_2$ with $B'$ via $f:S\to B'$. Set
$\shL_i=f_*(\shI_{C_i}/\shI_{C_i}^2)$. Then $\deg(\shL_i)= -C_i^2$
and $\shL_1\otimes\shL_2=\O_{B'}$. On the other hand, $\shT^1_X$ is
the dual of $\shL_1\otimes\psi^*(\shL_2)$. The kernel of the
homomorphism
$$
\phi_\shL : B'\lra \Pic^0(B'), \quad b\longmapsto
T^*_b(\shL)\otimes\shL^{-1}
$$
is the subgroup $B'_e\subset B'$ of $e$-torsion points
(\cite{Lange; Birkenhake 1992}, Lem.~4.7). Choose an origin $0\in
B'$ and some $b\in B'$ such that $\psi$ is the translation $T_b:B'\ra
B'$. So $\shT^1_X$ is trivial if and only if $b\in B'_e$, which means
$w|e$. \qed
\begin{remark}
\label{classify elliptic}
Here we see that a complete set of invariants of admissible $X$ with
nonrational algebraic normalization and positive degree is: The
$j$-invariant of $B'$, the degree $e>0$, and a translation
$\psi$ of $B'$ of finite order. In fact, $\PP^1$ bundles over $B'$
of degree $e>0$ with two disjoint sections are of the form
$\PP(\O_{B'}\oplus L^\vee)$ with $L$ a line bundle of degree $e$. Up to
pull-back by translations the latter are all isomorphic.
\end{remark}
We come to the construction of smoothings:
\begin{theorem}
\label{smoothing elliptic}
Suppose $X$ is an admissible surface (Definition~\ref{admissible})
with algebraic, nonrational normalization, of degree $e$ (possibly
$0$) and warp $w$. Then $X$ deforms to an elliptic principal bundle
of degree $e/w$ over an elliptic curve.
\end{theorem}
\proof In order to construct the desired smoothing, we first pass to a
Galois covering of $X$. The ruling $f:S\ra B'$ yields a bundle $X\ra
B$ over the isogenous elliptic curve $B=B'/\langle\psi\rangle$ whose
fibres are N\'eron $w$-gons. Let $G\subset\Pic^0(B')\subset\Aut(B')$
be the group of order $w$ generated by $\psi$. Consider the surface
$S'=S\times G$ and the isomorphisms
$$
\psi_j:C_2\times \left\{ \psi^j \right\}\lra C_1\times \left\{
  \psi^{j+1} \right\},\quad (s,\psi^j)\longmapsto (\psi(s),
\psi^{j+1})
$$
for $j\in\ZZ/w\/\ZZ$. The corresponding relation defines a normal
crossing surface $X'$ with $w$ irreducible components. Clearly, $X'$
is $d$-semistable with $K_{X'}=0$. The surface $S'$ is endowed with the
free $G$-action $\psi:(s,\psi^j)\mapsto(s,\psi^{j+1})$. It descends
to a free $G$-action on $X'$ with quotient $X=X'/G$.

We proceed similarly as in Theorem~\ref{smoothing Hopf}. Let $V$ be
the smooth toric variety corresponding to the infinite fan in $\ZZ^3$
generated by the cones
$$
\sigma_n= \Big\langle(0,n,1),(0,n+1,1)\Big\rangle,\quad n\in\ZZ.
$$
The projection $\pr_{3}:\ZZ^3\ra\ZZ$ defines a toric morphism
$V\ra\CC$. The special fibre $V_0$ is the product of a N\'eron
$\infty$-gon with $\CC^*$. The automorphism
$$
\Psi=
\begin{pmatrix}
  1&0&0 \\
  0&1&0\\
  0&1&1
\end{pmatrix}
\in \SL_3(\ZZ).
$$
acting from the right on row vectors maps $\sigma_n$ to
$\sigma_{n+1}$. By abuse of notation we also write $\Psi$ for the
induced automorphism of $V$. Choose $\alpha\in\Delta^*$ with
$B'=\CC^*/\langle \alpha^w\rangle$ and $B=\CC^*/\langle\alpha\rangle$.
The automorphism $(t_1,t_2,t_3)\mapsto (\alpha t_1,t_1^{e/w}t_2,t_3)$
of $(\CC^*)^3$ extends to another automorphism $\Phi$ of $V$. The
action of $\langle \Phi,\Psi\rangle$ on the preimage $U\subset V$ of
$\Delta\subset\CC$ is proper and free. Hence $\foX' = U/\langle
\Phi^w,\Psi^w\rangle$ is a smooth 3-fold, endowed with a projection
$\foX'\ra \Delta$. The general fibres $\foX'_t$, $t\in\Delta^*$, are
smooth Kodaira surfaces of degree $e$. The special fibre is
$\foX'_0=X'$. The action of $(\Phi\Psi)^w$ descends to a free
$G$-action on $\foX'$, which coincides with the given action on
$\foX'_0=X'$. Hence $\foX=\foX'/G$ is the desired smoothing. \qed

\medskip We head for the calculation of the versal deformation of $X$.
\begin{proposition}
\label{cohomology elliptic}
Suppose $K_X=0$. Then the following holds:
\renewcommand{\labelenumi}{(\roman{enumi})}
\begin{enumerate}
\item If $e>0$, then $h^0(\Theta_X)=1$, $h^1(\Theta_X)=2$, and
  $h^2(\Theta_X)=1$.
\item If $e=0$, then $h^0(\Theta_X)=2$, $h^1(\Theta_X)=3$, and
  $h^2(\Theta_X)=1$.
\end{enumerate}
\end{proposition}
\proof As in the proof of Lemma~\ref{cohomology Hopf} we use
$H^p( \Theta_X)\simeq H^{2-p}( \Omega_X^1/\tau_X) $ and the exact
sequence
$$
0 \lra \Omega_X^1/\tau_X \lra \Omega^1_S \lra \Omega^1_D \lra 0.
$$
Recall $h^{1,0}(S)=h^{1,0}(D)=1$. The map on the right in
$$
0\lra H^0(X,\Omega_X^1/\tau_X) \lra H^0(S,\Omega^1_S ) \lra
H^0(D,\Omega^1_D)
$$
is zero, since it is an alternating sum, so $h^2(\Theta_X)=1$.
Next we consider
$$
0 \lra H^0(D,\Omega^1_D) \lra H^1(X,\Omega_X^1/\tau_X) \lra
H^1(S,\Omega^1_S ) \lra H^1(D,\Omega^1_D).
$$
We have $h^{1,1}(S)=2 $ and $h^{1,1}(D)=1 $. The class of a fibre
$F\subset S$ in $H^{1,1}(S)$ maps to zero in $H^{1,1}(D)$. For $e=0
$, this also holds for the class of the section $C_1\subset S $, and
$h^1(\Theta_X)=3 $ follows. For $e>0 $, the image of the section does
not vanish, and $h^1(\Theta_X)=2 $ holds instead. Finally, the
sequence
$$
H^1(S,\Omega^1_S ) \lra H^1(D,\Omega^1_D) \lra
H^2(X,\Omega_X^1/\tau_X) \lra H^2(S,\Omega^1_S ) \lra 0
$$
is exact. Now $h^{1,2}(S)=1 $ yields $h^0(\Theta_X)= 2$ for $e=0
$, and $h^0(\Theta_X)= 1 $ for $e>0$. \endproof
\begin{corollary}
\label{tori elliptic}
Suppose $K_X=0$ and $e=0$. Then $X$ does not deform to a smooth
Kodaira surface.
\end{corollary}
\proof Write $S=\PP(\O_{B'}\oplus\shL)$ for some invertible
$\O_{B'}$-module $\shL$ of degree 0. Moving the gluing parameter
$\psi$ and the isomorphism class of $B'$ and $\shL$ gives a
3-dimensional locally trivial deformation of $X$. Since
$h^1(\Theta_X)=3$, the space parameterizing the locally trivial
deformations in the semiuniversal deformation $p:\foX\ra V$ of $X$ is
a smooth 3-fold. Now each fibre $\foX_t$ deforms to an elliptic
principal bundle of degree 0 (Thm.~\ref{smoothing elliptic}), and
the embedding dimension of $V$ is $\dim T^1_X=4$. This shows that $V$ is a
smooth 4-fold with an open dense set parameterizing complex tori.
Hence no fibre $\foX_t$ is isomorphic to a smooth Kodaira surface.
\qed
\begin{proposition}
\label{T^i elliptic}
Let $X$ be as in Theorem~\ref{smoothing elliptic} with $e>0$. Then $\dim
T^0_X=1$, $\dim T^1_X=3$ and $\dim T^2_X=2$.
\end{proposition}
\proof  In view of Proposition~\ref{cohomology elliptic} we only have
to show degeneracy of the spectral sequence of tangent cohomology
$E_2^{p,q}= H^p(\shT^q_X) \Rightarrow T^{p+q}_X$ at $E_2$ level, cf.\
Proposition~\ref{T^i Hopf}. This is the case iff $T^1_X\to
H^0(\shT^1_X)$ is surjective. Now by $d$-semistability
$H^0(\shT^1_X)$ is one-dimensional and any generator has no zeros.
Geometrically degeneracy of the spectral sequence therefore means the
existence of a deformation of $X$ over $\Delta_\varepsilon:=
\Spec\CC[\varepsilon]/\varepsilon^2$ that is not locally trivial.
This is what we established by explicit construction in
Theorem~\ref{smoothing elliptic}. More precisely, let $\foX \to
\Delta$ be a deformation of $X$. For $P\in D$ the image of the
Kodaira-Spencer class $\kappa\in T^1_X$ of this deformation in
$\shT^1_{X,P} \simeq \shExt^1 (\Omega^1_{X,P},\O_{X,P})$ is the
Kodaira-Spencer class of the induced deformation of the germ of $X$
along $D$. In appropriate local coordinates such deformations have
the form $xy-\varepsilon f(z)=0$ with $f$ inducing the section of
$\shExt^1 (\Omega^1_{X,P},\O_{X,P})$. Therefore $f\not\equiv 0$ iff
the total space $\foX$ is smooth at $P$. \qed
\begin{theorem}
\label{versal elliptic}
Suppose $X $ is admissible with algebraic, nonrational normalization,
of degree $e>0$ and warp $w$. Let $p:\shX\ra V $ be the semiuniversal
deformation of $X $. Then $V=V_1\cup V_2 $ has two irreducible
components, and the following holds:
\renewcommand{\labelenumi}{(\roman{enumi})}
\begin{enumerate}
\item $V_2 $ is a smooth surface and parameterizes the locally trivial
  deformations.
\item $V_1\cap V_2 $ is a smooth curve and parameterizes the
  $d$-semistable locally trivial deformations.
\item $V_1 $ is a smooth surface, and $V_1\setminus V_1\cap V_2 $
  parameterizes smooth Kodaira surfaces of degree $e/w$.
\end{enumerate}
\end{theorem}
\proof A similar situation has been found by Friedman in his study of
deformations of $d$-semistable K3 surfaces. We follow the proof of
\cite{Friedman 1983}, Theorem~5.10 closely.

Deformation theory provides a holomorphic map $h:T^1_X\ra T^2_X $ with
$V=h^{-1}(0) $, whose linear term is zero, and whose quadratic term is
given by the Lie bracket $1/2[v,v] $ (cf.\ e.g.\ \cite{Palamodov
  1976}). Moving the two parameters $j(E)\in \CC $, $\psi\in B' $ and
using $h^1(\Theta_X)=2 $, we see that the locally trivial deformations
are parameterized by a smooth surface $V_2\subset T^1_X $. Let
$h_2:T^1_X\ra \CC $ be a holomorphic map with $V_2=h_2^{-1}(0) $, and
$h_1:T^1_X\ra T_X^2 $ a holomorphic map with $h=h_1 h_2 $. Let
$L_1:T_X^1\ra T_X^2 $ and $L_2:T_X^1\ra \CC $ be the corresponding
tangential maps and set $V_1=h_1^{-1}(0)$.

We proceed by showing that $V_1\cap V_2 $ is a smooth curve, or
equivalently that the intersection of $\ker(L_1)$ with $H^1(\Theta_X)
= \ker(L_2) $ is 1-dimensional. The smoothing of $X $ constructed in
Theorem~\ref{smoothing elliptic} obviously has a smooth total space.
As in the proof of Proposition~\ref{T^i Hopf} we see that its
Kodaira-Spencer class $k\in T^1_X $ generates $\shT^1_X $.

Let $L:H^1(\Theta_X)\ra H^1(\shT_X^1)$ be the linear map $v\mapsto
[v,k] $. By \cite{Friedman 1983}, Proposition 4.5, its kernel are the
first order deformations which remain $d$-semistable. Since we can
destroy $d$-semistability by moving the gluing $\psi\in\Pic^0(E)$, this
kernel is 1-dimensional. As in \cite{Friedman 1983}, page 109, one
shows that for $v\in H^1(\Theta_X)$:
$$
L(v) = [v,k] = 1/2[v+k,v+k] = L_1(v+k)\cdot L_2(v+k) = L_1(v)\cdot
L_2(k).
$$
Using $L_2( k)\neq 0 $ we infer $ \ker L_1\cap H^1(\Theta_X)= \ker
L\simeq\CC$. It follows that $V_1\cap V_2$ is a smooth curve, which
parameterizes the $d$-semistable locally trivial deformations, and $V_1
$ must be a smooth surface. Our local computation together with the
interpretation of $j(E)\in\CC$ as coordinate on $V_1\cap V_2$ shows
that $V_1$ is nothing but the base of the smoothing of $X$ constructed
in Theorem~\ref{smoothing elliptic}. In particular, $V_1\setminus
V_1\cap V_2 $ parameterizes smooth Kodaira surfaces. \endproof

\section{Rational surfaces and honeycomb degenerations}

The goal of this section is to analyze rational admissible surfaces.
According to Proposition \ref{rational}, the ramification curve
$C\subset S $ of the normalization $\nu:S\ra X $ is a N\'eron 6-gon
$C=C_0\cup\ldots \cup C_5 $ with $K_S=-C$. We suppose $C_i\cap
C_{i+1}\neq\emptyset$, regarding the indices as elements in
$\ZZ/6\ZZ$. The singular locus $D\subset X$ is the seminormal curve
with normalization $\PP^1\times\ZZ/3\ZZ$ and imposed relations
$(0,i)\sim(0,j)$ and $(\infty,i)\sim(\infty,j)$ for $i,j\in\ZZ/3\ZZ$.
The gluing $\varphi:C\ra D $ identifies pairs of irreducible
components in $C $. Since the resulting surface should be normal
crossing the fibers of $S\to X$ have cardinality at most $3$. There
are two possibilities left. The first alternative is
$\varphi(C_i)=\varphi(C_{i+3})$ and $\varphi(C_i\cap
C_{i+1})=\varphi(C_{i+2}\cap C_{i+3})$ for all $i\in\ZZ/6\ZZ $. We
call it \textit{untwisted} gluing. The second alternative is
$\varphi(C_i)=\varphi(C_{i+3})$, $\varphi(C_{i+1})=\varphi(C_{i-1})$,
$\varphi(C_{i+2})=\varphi(C_{i-2})$ and $\varphi(C_i\cap
C_{i+1})=\varphi(C_{i-1}\cap C_{i-2})$ for some $i\in\ZZ/6\ZZ $. Call
it \emph{twisted} gluing. The distinction is important for the
canonical class:
\begin{proposition}
\label{untwisted rational}
The condition $K_X=0$ holds if and only if the gluing $\varphi:C\ra D$
is untwisted.
\end{proposition}
\proof As in Propositions~\ref{untwisted Hopf} and \ref{untwisted
  elliptic} the condition $K_X=0$ is equivalent to $H^2(\O_{X})\neq
0$. By the exact sequence
$$
H^1(S,\O_{S})\oplus H^1(D,\O_{D}) \lra H^1(C,\O_{C}) \lra
H^2(X,\O_{X}) \lra 0
$$
and with $h^{0,1}(S)=0$ this holds precisely if
$\varphi^*:H^1(\O_{D}) \ra H^1(\O_{C})$ vanishes. The bicoloured
graphs attached to the curves $C,D$ (\cite{Deligne; Rapoport 1972},
section 3.5) are
\begin{center}
\setlength{\unitlength}{0.0004in}%
\begingroup\makeatletter\ifx\SetFigFont\undefined%
\gdef\SetFigFont#1#2#3#4#5{%
     \reset@font\fontsize{#1}{#2pt}%
     \fontfamily{#3}\fontseries{#4}\fontshape{#5}%
     \selectfont}%
\fi\endgroup%
\begin{picture}(5566,3119)(1418,-3200) \thicklines
\put(1501,-961){\circle{150}} \put(1501,-1861){\circle{150}}
\put(2401,-2461){\circle{150}} \put(3301,-1861){\circle{150}}
\put(3301,-961){\circle{150}} \put(3301,-1411){\circle*{150}}
\put(1501,-1411){\circle*{150}} \put(2851,-661){\circle*{150}}
\put(1951,-661){\circle*{150}} \put(1951,-2161){\circle*{150}}
\put(2851,-2161){\circle*{150}} \put(6001,-361){\circle{150}}
\put(6001,-1411){\circle{150}} \put(6001,-2461){\circle{150}}
\put(5101,-1411){\circle*{150}} \put(2401,-361){\circle{150}}
\put(6901,-1411){\circle*{150}}
\put(5700,-3100){\makebox(0,0)[lb]{$\Gamma(D)$}}
\put(1576,-886){\line( 5, 3){750}}
\put(2476,-436){\line( 5,-3){750}} \put(1501,-1036){\line( 0,-1){750}}
\put(1576,-1936){\line( 5,-3){750}} \put(2476,-2386){\line( 5,
3){750}} \put(3301,-1786){\line( 0, 1){750}} \put(5101,-1411){\line(
5, 6){817.623}} \put(5176,-1411){\line( 1, 0){750}}
\put(5176,-1486){\line( 5,-6){750}} \put(6076,-436){\line(
5,-6){817.623}} \put(6076,-1411){\line( 1, 0){825}}
\put(6076,-2386){\line( 5, 6){817.623}} \put(3901,-1411){\vector( 1,
0){600}} \put(2100,-3100){\makebox(0,0)[lb]{ $\Gamma(C)$}}
\end{picture}
\end{center}
Here the white vertices are the irreducible components, the black
vertices are the singularities, and the edges are the ramification
points on the normalization. Since the curves in question are
seminormal, the map $\varphi^*:H^1(\O_{D}) \ra H^1(\O_{C})$ coincides
with $\varphi^*:H^1(\Gamma(D),\CC)\ra H^1(\Gamma(C),\CC)$. A direct
calculation left to the reader shows that it vanishes if
and only if $\varphi$ is untwisted. \qed

\medskip From now on we assume that the gluing $\varphi:C\ra D$ is
untwisted. According to Proposition \ref{rational}, the surface $S$
is obtained from a Hirzebruch surface $f':S'\ra\PP^1$ by blowing-up
twice. Our conventions are that $C_0,C_3\subset S$ are sections of
the induced fibration $f:S\ra\PP^1$, that $C_1\cup C_2$, $C_4\cup C_5$
are the fibres over $0,\infty\in\PP^1$, and that $C_1,C_4$ are
exceptional for the contraction $r:S\ra S'$. Let us also assume that
the image of $C_0$ in $ S'$ is a minimal section. If $e\geq 0$ is the
degree of the Hirzebruch surface $S'$, this means $C_0^2=-e-1$ and
$C_3^2=e-1$.

Set $D_i=\varphi(C_i)$ for $i\in\ZZ/3\ZZ$. The space of all untwisted
gluings $\varphi:C\ra D$ is a torsor under
\begin{equation}
\label{gluing parameter rational}
\Aut^0(D)=\prod\Aut^0(D_i) \simeq (\CC^*)^3.
\end{equation}
The action of $\phi\in\Aut^0(D)$ is given by composing $\varphi\mid
C_0\cup C_1\cup C_2$ with $\phi$. We call $\varphi\mid C_0\cup C_3$
the \emph{vertical} gluing. The ruling $f:S\ra \PP^1$ yields a
preferred vertical gluing, which identifies points of the same fibre.
Every other vertical gluing differs from the preferred one by a
\textit{vertical gluing parameter} $\zeta\in\CC^*$. We call its order
$w\geq0$ the \emph{warp} of $X$. The warp is important for the
calculation of $\shT^1_X$:
\begin{proposition}
\label{d-semistable rational}
The surface $X$ is $d$-semistable if and only if the warp $w$ divides
the degree $e$.
\end{proposition}
\proof The inclusions $D_i\cup D_j\subset D$ give an injection $
\Pic(D)\subset \prod_{i\neq j} \Pic(D_i\cup D_j)$. We proceed by
calculating the class of $\shT^1_X\mid D_0\cup D_1$. Consider the
normal crossing surface $\bar{S}$ defined by the cocartesian diagram
$$
\begin{CD}
  C_2\cup C_5 @>>> S \\
  @VVV @VVV\\
  D_2 @>>> \bar{S }
\end{CD}
$$
Let $\bar{C}_i\subset \bar S$ be the images of $C_i\subset S$. The
ideals $\shI_{01},\shI_{34}\subset\O_{\bar{S}}$ of the Weil divisors
$\bar{C}_{01} = \bar{C}_0\cup \bar{C}_1$ and $\bar{C}_{34}=
\bar{C}_3\cup \bar{C}_4$ are invertible. A local computation shows
that
$$
(\shT^1_X \mid D_0 \cup D_1)^\vee \simeq \shI_{01}/\shI_{01}^2
\otimes \psi^*(\shI_{34}/\shI_{34}^2),
$$
where $\psi:\bar{C}_{01}\lra \bar{C}_{34}$ is the isomorphism
obtained from the induced gluing $\overline{\varphi}:
\bar{C}_{01}\cup\bar{C}_{34}\ra D_0\cup D_1$. Another local
computation gives
$$
\shI_{01}/\shI_{01}^2 \mid \bar{C}_{0} \simeq
\shI_{C_0}/\shI_{C_0}^2(-C_0\cap C_{-1}), \quad \shI_{01}/\shI_{01}^2
\mid \bar{C}_{1} \simeq \shI_{C_1}/\shI_{C_1}^2(-C_0\cap C_{1}).
$$
It follows that $\shT^1_X\mid D_0\cup D_1$ is trivial if and only
if $\zeta^e=1$. A similar argument applies to $D_0\cup D_2$ and
$D_1\cup D_2$. \qed
\begin{remark}
\label{classify rational}
The previous two propositions yield the following classification
of admissible $X$ with rational normalization and positive degree
$e$: According to Equation~\ref{gluing parameter rational},
the isomorphism class of $X$ is determined by $e>0$ and 3 gluing parameters.
By Proposition~\ref{d-semistable rational} the vertical one is an $e$-th
root of unity $\zeta$. For $(e,\zeta)$ fixed the automorphisms of $S$ act
on the remaining two horizontal gluing parameters with quotient
isomorphic to $\CC^*$.
\end{remark}
We come to the existence of smoothings:
\begin{theorem}
\label{smoothing rational}
Suppose $X$ is an admissible surface (Definition~\ref{admissible})
with rational normalization, of degree $e$ (possibly $0$) and warp
$w$. Then $X$ deforms to an elliptic principal bundle of degree $e/w$
over an elliptic curve.
\end{theorem}
\proof It should not be too surprising that the construction is a
modification of both the constructions in Theorems~\ref{smoothing
Hopf} and \ref{smoothing elliptic}.

First, we simplify matters by passing to a Galois covering. Let
$G\subset \CC^*$ be the group of order $w$ generated by the vertical
gluing parameter $\zeta\in\CC^*$. Since $w| e$, the $G$-action
$z\mapsto \zeta z$ on $\PP^1$ lifts to a $G$-action on the Hirzebruch
surface $S'$, and hence to the blowing-up $S$. The diagonal action on
$\tilde{S}=S\times G$ is free. Let $\psi_i:C_i\ra C_{i+3}$ be the
isomorphisms induced by the gluing $\varphi:C\ra D$. This gives
$G$-equivariant isomorphisms
$$
\psi_{ij}: C_i\times \left\{ \zeta^j \right\} \lra C_{i+3}\times
\left\{\zeta^ {j+1} \right\}, \quad (s,\zeta^j)\longmapsto (\psi_i(s),
\zeta^{j+1}).
$$
Let $X'$ be the normal crossing surface obtained from $\tilde S$
modulo the relation imposed by the $\psi_{ij}$. Then $G$ acts freely
on $X'$ with quotient $X=X'/G$.

The next task is to construct a $G$-equivariant smoothing of $X'$.
Again toric geometry enters the scene. Consider the infinite fan in
$\ZZ^4$ generated by the cones
\begin{eqnarray*}
\sigma_{m,n} &=& \Big\langle(m,e\cdot{\textstyle \binom{m}{2}},1,0),(m+1,e
\cdot{\textstyle \binom{m+1}{2}},1,0), (0,n,0,1),(0,n+1,0,1)\Big\rangle
\end{eqnarray*}
for $m,n\in\ZZ$. Let $V$ be the corresponding 4-dimensional smooth
toric variety. The projection $\pr_{34}:\ZZ^4\ra\ZZ^2$ defines a
toric morphism $\pr_{34}:V\ra \CC^2 $. The special fibre $V_0$ is
isomorphic to a bundle of N\'eron $\infty$-gons over the N\'eron
$\infty$-gon. Let $f:\hat{\CC}^2\ra\CC^2$ be the blowing-up of the
origin $0\in\CC^2$ with the closed toric orbits removed and $E \subset
\hat{\CC}^2$ the exceptional set. It is isomorphic to $\CC^*$. The
cartesian diagram
$$
\begin{CD}
  \hat{ V } @>>> V \\
  @VVV @VVV\\
  \hat{\CC}^2 @>>> \CC^2
\end{CD}
$$
defines a smooth toric 4-fold $\hat{ V}$, which is an open subset
of the blowing-up of $ V$ with center $ V_0\subset V$. The
exceptional divisor $\hat{ V}_E=\hat{ V}_0$ of $\hat{ V}\ra V$ is
isomorphic to $E\times V_0 $. The exceptional divisor $E$ will be a
parameter space of the whole construction via the $\CC^*$-bundle
$\hat{\CC}\ra E$ induced by the canonical map $\CC^2\setminus 0\ra
\PP^1$.

Let $Z\subset \hat{ V}$ be the set of 1-dimensional toric orbits.
Each fibre $Z_t$, $t\in E$, consists of the discrete set of
non-normal-crossing singularities in $\hat{ V}_t \simeq V_0$, as
illustrated by Figure~1. Let $\bar{ V}\ra\hat{ V}$ be the blowing-up
with center $Z\subset \hat{ V}$. The exceptional divisor
$\bar{Z}\subset\bar{ V}$ is easy to determine. Each fibre $\hat V_t$,
$ t\in E$, is smooth, except for the points in $Z_t$. At such points
a local equation for $\hat V_t$ is $T_1T_2=tT_3T_4$, which is an
affine cone over a smooth quadric in $\PP^2$. Hence $\bar{Z}$ is a
disjoint union of smooth quadrics in $ \PP^3$, each isomorphic to
$\PP^1\times\PP^1$. The whole exceptional divisor is $\bar{Z}\simeq
E\times Z\times \PP^1\times\PP^1$. The picture over a fixed $t\in E$
is
depicted in Figure~2.\\[3ex]
\mbox{}\hfill\includegraphics{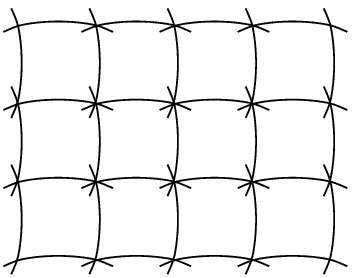}\hfill
\includegraphics{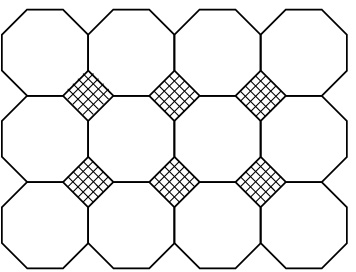}\hfill\mbox{}\\[2ex]
\mbox{}\hspace{1.7cm} Fig.\ 1\hspace{4.2cm} Fig.\ 2\\[3ex]
Here the octagons lie on the strict transform of $ V_t$, and the
squares are contained in the exceptional divisor $\bar{Z}_t$.

We seek to contract $\bar{Z}\subset \bar{ V}$ along one of the two
rulings of $\PP^1\times\PP^1$. Let $F\subset \bar{Z}$ be a
$\PP^1$-fibre. A local calculation shows that $\O_F(\bar Z)$ has
degree $-1$, and that the Cartier divisors $K_{\bar V}$ and $\O_{\bar
  V}(\bar Z)$ coincide in a neighborhood of $F\subset\bar{ V}$. So
the Nakano contraction criterion \cite{Nakano 1971} applies, and there
exists a contraction $\bar{ V}\ra\tilde{ V}$ which restricts to the
projection
$$
\pr_{123}: \bar{Z}=E\times Z\times \PP^1\times\PP^1 \lra E\times
Z\times \PP^1.
$$
The special fibre $\tilde{ V}_0$ resembles an infinite system of
honeycombs:\\[3ex]
\centerline{\includegraphics{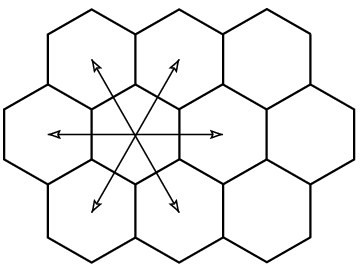}}\\[-2ex]
\mbox{}\hspace{2.5cm} Fig.\ 3\\[3ex]
Finally, we want to arrive at compact surfaces, so we seek to divide
out a cocompact group action. Consider the two commuting
automorphisms
$$
\Phi=
\begin{pmatrix}
  1&e&0&0 &0\\
  0&1&0&0&0 \\
  1&0&1&0&0\\
  0&0&0&1 &0\\
  0&1&0&0&1 \\
\end{pmatrix}
\quadand \Psi=
\begin{pmatrix}
  1&0&0&0&0\\
  0&1&0&0&0\\
  0&0&1&0&0\\
  0&1&0&1&0\\
  0&0&0&0&1\\
\end{pmatrix}
$$
of $\ZZ^5$ acting from the right on row vectors. We regard our fan
in $\ZZ^4$ as a fan in $\ZZ^5=\ZZ^4\oplus\ZZ$. Torically, the trivial
$\ZZ$-factor amounts to going over to $V\times\CC^*$, which will give
a horizontal gluing parameter. Since $\Phi(\sigma_{m,n})
=\sigma_{m+1,n}$ and $\Psi(\sigma_{m,n}) =\sigma_{m,n+1}$, we get an
induced action of the discrete group $\ZZ^2=\langle \Phi,\Psi\rangle$
on $V\times\CC^*$. Let $U\subset V$ be the preimage of $\Delta^2
\subset\CC^2$. The action is proper and free on $U\times\CC^*$. Let
$\tilde{U} \subset\tilde{V}$ be the corresponding preimage. It is
easy to see that the action on $V\times \CC$ induces an action on
$\tilde{ V}\times\CC^*$ which is proper and free on $\tilde{
  U}\times\CC^*$. So the quotient $ \foX'= (\tilde{U}\times\CC^*)/
\langle \Phi,\Psi^w\rangle$ is a smooth complex 5-fold. The action on
the special fibre is indicated in Figure~3 by the arrows. The general
fibres $ \foX'_t$, $(t,\lambda)\in f^{-1}(\Delta^*\times\Delta^*)
\times\CC^*$ are elliptic bundles of degree $e$ over the elliptic
curve $\CC^*/\langle t_1\rangle$ with fibre $\CC^*/\langle
t_2\rangle$, where $f(t)=(t_1,t_2)$. Some special fibre $\foX'_t$,
$t\in E\times\CC^*$, is isomorphic to $X'$, since $t_5$ moves through
all horizontal gluings.

It remains to extend the $G$-action on $X'$ to a $G$-action on
$\foX'$. The automorphism
$$
(\CC^*)^5\lra(\CC^*)^5,\quad (t_1,t_2,t_3,t_4,t_5) \longmapsto
(\zeta t_1,t_2t_4,t_3,t_4,t_5)
$$
of the torus extends to an automorphism of the torus embedding $V$.
It induces a free $G$-action on $\foX'$, which is the desired
extension. Set $\foX=\foX'/G$. General fibres $ {\foX}_t$,
$t\in(\CC^2\setminus0)\times\CC^*$ are elliptic bundles of degree
$e/w$ over an elliptic curve. Some special fibre ${\foX} _t$, $t\in
E\times\CC^*$ is isomorphic to $X$. \qed

\medskip We turn our attention to the versal deformation of $X$ and
calculate the relevant cohomology groups:
\begin{proposition}
\label{cohomology rational}
Suppose $K_X=0$. Then the following holds:
\renewcommand{\labelenumi}{(\roman{enumi})}
\begin{enumerate}
\item If $e>0$, then $h^0(\Theta_X)=1$, $h^1(\Theta_X)=2$, and
  $h^2(\Theta_X)=0$.
\item If $e=0$, then $h^0(\Theta_X)=2$, $h^1(\Theta_X)=3$, and
  $h^2(\Theta_X)=0$.
\end{enumerate}
\end{proposition}
\proof As in Lemmas~\ref{cohomology Hopf} and \ref{cohomology
  elliptic} we use $H^p(X,\Theta_X) \simeq
H^{2-p}(X,\Omega^1_X/\tau_X)$ and the exact sequence
$$
0 \lra \Omega^1_X/\tau_X \lra \Omega_S^1 \lra \Omega^1_{\tilde{D}}
\lra 0.
$$
Here $\tilde{D}$ is the normalization of $D$. The inclusion $
H^0(\Omega^1_X/\tau_X) \subset H^0(\Omega_S^1)$ and $b_1(S)=0$ gives
$h^2(\Theta_X)=0$. Next consider the exact sequence
$$
0 \lra H^1(X,\Omega^1_X/\tau_X) \lra H^1(S,\Omega^1_S)\lra
H^1(\tilde{D},\Omega^1_{\tilde{D}})
$$
The task is to determine the rank of the map on the right. This is
done as in \cite{Friedman 1983}, p.~91: One has to compute the images
in $H^{1,1}(\tilde{D})$ of the classes of $C_i\subset S$ using
intersection numbers on $S$. For $e>0$ this gives $h^1(\Theta_X)=2$,
whereas for $e=0$ the result is $h^1(\Theta_X)=3$. We leave the
actual computation to the reader. Finally, the exact sequence
$$
H^1(S,\Omega^1_S)\lra H^1(\tilde{D},\Omega^1_{\tilde{D}}) \lra
H^2(X,\Omega^1_X/\tau_X) \lra H^2(S,\Omega^1_S),
$$
together with the preceding observations and $b_3(S)=0$ gives the
stated values for $h^0(\Theta_X)$. \qed
\begin{corollary}
\label{tori rational}
Suppose $K_X=0$ and $e=0$. Then $X$ does not deform to a smooth
Kodaira surface.
\end{corollary}
\proof The locally trivial deformations of $X$ are unobstructed since
$h^2(\Theta_X)=0$. They define a smooth 3-fold $V'\subset V$ in the
base of the semiuniversal deformation $\foX\ra V$. It is easy to see
that it is given by the three gluing parameters in $\varphi:C\ra D$.
Since each fibre $\foX_t$ deforms to a complex torus
(Thm.~\ref{smoothing rational}), $V$ is smooth of dimension $4$ and
no fibre can be isomorphic to a smooth Kodaira surface. \qed
\begin{proposition}
\label{T^i rational}
Let $X$ be as in Theorem~\ref{smoothing rational} with $e>0$. Then $\dim
T^0_X=1$, $\dim T^1_X=3$ and $\dim T^2_X=2$.
\end{proposition}
\proof
This is shown as the similar statement of Proposition~\ref{T^i elliptic}.
\qed
\begin{theorem}
\label{versal rational}
Suppose $X$ is admissible with rational normalization, of degree $e>0$
and warp $w$. Let $p:\foX\ra V$ be the semiuniversal deformation of
$X$. Then $V=V_1\cup V_2$ consists of two irreducible components, and
the following holds: \renewcommand{\labelenumi}{(\roman{enumi})}
\begin{enumerate}
\item $V_2$ is a smooth surface parameterizing the locally trivial
  deformations.
\item $V_1\cap V_2$ is a smooth curve parameterizing the locally
  trivial deformations which remain $d$-semistable.
\item $V_1\setminus V_1\cap V_2$ parameterizes smooth Kodaira surfaces
  of degree $e/w$.
\end{enumerate}
\end{theorem}
\proof The argument is as in the proof of Theorem~\ref{versal
  elliptic}. \qed

\section{The completed moduli space and its boundary}

In this final section, we take up a global point of view and analyze
degenerations of smooth Kodaira surfaces in terms of moduli spaces.
Here we use the word ``moduli space'' in the broadest sense, namely
as a topological space whose points correspond to Kodaira surfaces.
We do not discuss wether it underlies a coarse moduli space
or even an analytic stack or analytic orbispace.

Let $\foK_d=\CC\times\Delta^*$ be the moduli space of smooth Kodaira
surfaces of degree $d>0$, and $\foK=\cup_{d>0}\foK_d$ their union.
The points of $\foK$ correspond to the isomorphism classes $[X]$ of
smooth Kodaira surfaces. The topology on $\foK$ is induced by what we
suggest to call the \emph{versal topology} on the set $\foM$ of
isomorphism classes of all compact complex spaces. The versal
topology is the finest topology on $\foM$ rendering continuous all
maps $V\ra \foM$ defined by flat families $\foX\ra V$ which are versal
for all fibres. The complex structure on $\foK$ comes from Hodge
theory: One can view $\foK_d$ as the period domain of polarized pure
Hodge structures of weight 2 on $H^2(X,\ZZ)/\text{torsion}$ divided by
the automorphism group of this lattice. According to \cite{Borcea
  1984}, the induced structure of a ringed space on $\foK_d$ is the
usual complex structure on $\CC\times \Delta^*$.

Let $\foK\subset\overline{\foK}$ be the space obtained by adding all
admissible surfaces in the closure of $\foK\subset\foM$. The surfaces
parameterized by the boundary $\foB=\overline{\foK}\setminus\foK$ are
called \emph{d-semistable} Kodaira surfaces. These are nothing but
the admissible surfaces deforming to smooth Kodaira surfaces. The
boundary decomposes into three parts
$$
\foB = \foB^h\cup\foB^r\cup \foB^e
$$
according to the three types of admissible surfaces. Here $\foB^h$
refers to the surfaces whose normalization are Hopf surfaces, $\foB^r$
to the surfaces with rational normalization, and $\foB^e$ to surfaces
whose normalization is ruled over an elliptic base. We refer to
these parts of $\overline\foK$ as Hopf, rational and elliptic ruled
stratum respectively.
\begin{proposition}
\label{boundary components}
The irreducible components of the boundary $\foB$ are smooth complex
curves. The components in $\foB^h$ are isomorphic to
$\Delta^*\simeq\left\{ \infty \right\}\times\Delta^*$, the components
in $\foB^r$ are isomorphic to $\CC^*$, and the components in $\foB^e$
are isomorphic to $\CC\simeq\CC\times \left\{ 0 \right\}$.
\end{proposition}
\proof
This follows from Remarks~\ref{classify Hopf}, \ref{classify elliptic}
and \ref{classify rational}.
\qed

\medskip
Locally, the completion $\foK\subset\overline{\foK}$ is isomorphic to
the blowing-up of $(\infty,0)\in\PP^1\times\Delta$, but with the
points $0,\infty\in\PP^1$ on the exceptional divisor removed. This
follows in particular from the construction in Theorem~\ref{smoothing
rational} of a family over the blow up of $\CC^2$ with $2$ points
removed, with fiber over $(t_1,t_2)\in \CC^*\times\CC^*$ a smooth
Kodaira surface with invariants $(j,\alpha)= (\exp(t_1), t_2)$. To
complete further we need to enlarge the class of generalized Kodaira
surfaces. The least one would hope for is that any family $\foX^*\to
\Delta^*$ of generalized Kodaira surfaces, that can be completed to a
proper family over $\Delta$, can be completed by a generalized
Kodaira surface.

We first show that this is impossible if we consider only normal
crossing surfaces.

\begin{theorem}
Let $\foX\to \Delta$ be a degeneration of $d$-semistable Kodaira
surfaces with elliptically ruled normalization. Assume that $\foX$ is
bimeromorphic to a K\"ahler manifold, and that the $j$-invariant of
the base of $\foX_t$ tends to $0$ with $t\in\Delta$. Then $\foX_0$ is
not of normal crossing type.
\end{theorem}

\proof
Let $\tilde\foX\to \foX$ be the normalization. Let $\foB$ be the
component of the Douady space of holomorphic curves in $\tilde\foX$
that contain a fiber of $\tilde\foX_t\to B_t$ for general
$t\in\Delta$. Since these curves are contained in fibers of
$\tilde\foX\to \Delta$ there is a map $\foB\to\Delta$. By the
K\"ahler assumption this map is proper \cite{Fujiki 1984}. Let
$\tilde\foX'\to\foB$ be the universal family. The universal map
$\tilde\foX'\to\tilde \foX$ is an isomorphism over $\Delta^*$, hence
bimeromorphic. By desingularisation we may dominate this map by
successively blowing up $\tilde\foX$ in smooth points and curves.
This can be arranged to keep the property ``normal crossing''. In
other words, we can assume that $\tilde\foX'=\tilde\foX$ is a
fibration by rational curves over a degenerating family
$\foB\to\Delta$ of elliptic curves.

Let $\foC',\foC''\subset\tilde \foX$ be the two components of the
preimage of the singular locus of $\foX$ mapping onto $\Delta$. For
any $b\in \foB$ the corresponding rational curve $F_b$ intersects
$\foC', \foC''$ in one point each. Let $b_0\in\foB_0$ be a node of the
nodal elliptic curve over $0\in\Delta$. Then each intersection
$\foC'\cap F_{b_0}$ and $\foC''\cap F_{b_0}$ gives a point of
multiplicity at least $2$ on $\tilde\foX_0$. As $\foC',\foC''$ are
identified under $\tilde\foX\to\foX$ this leads to a point of
multiplicity at least $4$ on $\foX_0$.
\qed
\medskip

There are however various completions if we admit \emph{products of
normal crossing singularities}. In dimension $2$ the only such
singularity that is not normal crossing is a point $(X,x)$ of multiplicity
$4$. It has as completed local ring
$$
\O_{X,x}^\wedge = \CC[[T_1,\ldots,T_4]] /(T_1T_2,T_3T_4).
$$
In particular, it is a complete intersection and hence still
Gorenstein. We will call $(X,x)$ a \emph{quadrupel point}. The
singular locus $C\subset X$ consists of the four coordinate lines
$C_1,\ldots, C_4$, and $C_i, C_j$ lie on the same irreducible
component iff $i-j\not \equiv 2$ modulo $4$. The embedding dimension of
$(X,x)$ is $4$. It is therefore not embeddable into a smooth 3-fold.
Its appearance here is perhaps not so surprising, as it occurs also
in certain compactifications of moduli of polarized abelian surfaces
\cite{Nakamura 1975}.

We were nevertheless unable to select a natural class of generalized
Kodaira surfaces that would satisfy the mentioned completeness
property. Therefore we content ourselves to define two surfaces with
quadrupel points, connecting the elliptic ruled stratum to the Hopf
stratum and to the rational stratum respectively. A surface
connecting the Hopf stratum with the rational stratum could not
be found, although such surfaces should probably exist.

\medskip Fix an integer $e> 0$. Let $S$ be the Hirzebruch surface of
degree $e$. Choose two sections $C_0, C_2\subset S$ with $C_0^2=-e$,
$C_2^2=e$, and let $C_1, C_3$ be the fibers over $0,\infty\in \PP^1$.
We define the surface $X_1$ by gluing $C_1, C_3$ by any isomorphism
(all choices give isomorphic results), and $C_0, C_2$ by an
isomorphism of finite order $w$ over $\PP^1$. Note that $X_1$ is
normal crossing except at one quadrupel point. As in the proof of
Proposition~\ref{untwisted rational} one checks  $K_{X_1}=0$.

\begin{proposition}
Suppose $w$ divides $e$. There is a family $\foX\to \CC^2$ with
the following properties:
\begin{enumerate}
\item $\foX_0\simeq X_1$
\item the fibers over $\CC^*\times\{0\}$ are $d$-semistable
Kodaira surfaces of warp $w$ and degree $e$ with elliptic ruled
normalization
\item the fibers over $\{0\}\times\CC^*$ are
$d$-semistable Kodaira surfaces of warp $w$ and degree $e$ with
nonalgebraic normalization
\end{enumerate}
\end{proposition}

\proof In the proof of Theorem~\ref{smoothing rational} we
constructed a toric morphism $V\to\CC^2$ with central fiber a bundle
of N\'eron $\infty$-gons over the N\'eron $\infty$-gon. The
restriction of $\Phi,\Psi\in \GL(\ZZ^5)$ defined there to
$\ZZ^4\simeq\ZZ^4\oplus\{0\}$ defines an action of $\ZZ^2$ on the fan
defining $V$. The induced action on $V$ is proper and free, and
commutes with the map to $\CC^2$. We obtain a proper family $\foX'=
V/\ZZ^2\to \CC^2$. By construction of $V$ there is also an action of
$G\simeq \ZZ/w\ZZ$ on $\foX'/\CC^2$. The family $\foX'/G\to \CC^2$
has the desired properties.
\qed

\medskip Finally we study a surface connecting the rational and the
elliptic stratum, under the expense of changing the warp. Again we
fix an integer $e> 0$. Let $S=S'\cup S''$ be the disjoint union of
two Hirzebruch surfaces of degrees $e+1$ and $e-1$ respectively.
Choose two sections $C'_0,C_2'\subset S'$ with $(C'_0)^2=-(e+1)$ and
$(C'_2)^2=e+1$, and let $C_1',C_3'$ be the fibres over
$0,\infty\in\PP^1$. Make the analogous choices $C''_0,C_2''$ etc.\
for $S''$.  We glue $C_i'$ to $C_{i+2}''$ in a way that preserves the
intersection points indicated in the following figure by circles and
crosses.\\[4ex]
\centerline{\includegraphics{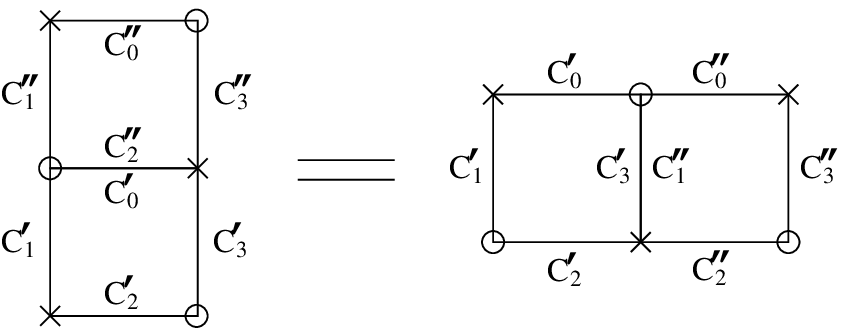}}\\[-2ex]
\mbox{}\hspace{0.5cm} Fig.\ 4\\[3ex]
For odd $i$ the isomorphism $C_i'\to C_i''$ can be chosen
arbitrarily; for even $i$ we glue compatibly with the projections
$S'\to \PP^1$, $S''\to \PP^1$. The result is a reduced surface $X_2$
that is normal crossing except at two points corresponding to the
circles and the crosses in the figure. It is a bundle of N\'eron
$1$-gons over $\PP^1$. Again as in the proof of Proposition
\ref{untwisted rational}, one verifies that $K_{X_2}=0$ holds.
\begin{proposition}
\label{deformation quadrupel}
\begin{enumerate}
\item $X_2$ deforms to $d$-semi\-stable Kodaira surfaces with
rational normalization of degree $e$ and warp $1$.
\item
$X_2$ deforms to $d$-semistable Kodaira surfaces
with elliptic ruled normalization of degree $2e$ and warp $2$.
\end{enumerate}
\end{proposition}
\proof (i) The idea is to modify a fiberwise degeneration of
Hirzebruch surfaces. Let $\widetilde S'$, $\widetilde S''$ be
Hirzebruch surfaces of degree $e$. Define curves $\widetilde{C}'_i
\subset \widetilde{S'}$ and $\widetilde{C}''_i \subset
\widetilde{S''}$ as above.

Let $p:\widetilde{\foZ}\ra\Delta$ be a flat family whose general
fibres $\widetilde{\foZ}_t$, $t\neq0$, are Hirzebruch surfaces of
degree $e$; the closed fiber $\widetilde\foZ_0$ is the union of
$\widetilde S'$ and $\widetilde S''$ with $\widetilde{C}_2'$ and
$\widetilde{C}_0''$ identified, analogously to the left-half of
Figure~4. Such a family can be constructed either by toric methods or
as follows. Let $L\to \PP^1$ be the $\CC^*$-bundle of degree $e$. Let
$\foC\to \Delta$ be a versal deformation of the nodal rational curve
$xy=0$ in $\PP^2$. The action $t\cdot (z,w)= (tz, t^{-1}w)$ of
$\CC^*$ on $\foC_0$ extends to an action on $\foC$ over $\Delta$. We
may take $\widetilde \foZ= (L\times \foC)/\CC^*$ with $\CC^*$ acting
diagonally.

The points $\widetilde{C}_0'\cap \widetilde{C}_1'$ and
$\widetilde{C}_2''\cap \widetilde{C}_3''$ on $\widetilde{\foZ}_0$ can
be lifted to sections of $\widetilde{\foZ}\ra\Delta$. Let
$\foZ\ra\widetilde{\foZ}$ be the blowing-up of these sections. The
strict transform of $\widetilde{C}_0'\cup \widetilde{C}_2''$ is a
Cartier divisor on $\foZ_0$. It can be extended to an effective
Cartier divisor on $\foZ$, whose associated line bundle we denote by
$\shL$. Consider the factorization $\foZ\to \Delta\times \PP^1$ of
$\foZ\to\Delta$. The relative base locus of $\shL$ over
$\Delta\times\PP^1$ is obviously finite. According to the
Fujita-Zariski Theorem (\cite{Fujita 1983}, Thm.\ 1.10) the
corresponding invertible $\O_{\foZ}$-module $\shL$ is relatively
semiample over $\PP^1\times\Delta$. Moreover, $\shL$ is relatively
ample over $\PP^1\times\Delta^*$. Let $\foZ\ra\foS$ be the
corresponding contraction. The exceptional locus of this contraction
is the strict transform of $\widetilde{C}_1'\cup \widetilde{C}_3''$.

Now $\foS\to \Delta$ has as general fiber the blowing up of a
Hirzebruch surface of degree $e$ in two points as needed for the
construction of $d$-semistable Kodaira surfaces with rational
normalization. The central fiber consists of a union of two
Hirzebruch surfaces of degrees $e+1$ and $e-1$ respectively, with a
section of degree $e+1$ glued to a section of degree $-(e-1)$. So
this is a partial normalization of $X_2$. The gluing morphism
$\varphi:C\ra D$ on the special fibre $\foS_0=S$ extends to a gluing
morphism over $\Delta$. The corresponding cocartesian diagram gives
the desired flat family $\foX\ra\Delta$ with $\foX_0=X_2$.

\noindent (ii) For the elliptic case consider the partial
normalization $Y$ of $X_2$ obtained from $S$ by identifying $C'_1$
with $C''_3$ and $C'_3$ with $C''_1$, see the right half of Figure~4.
It is the projective closure of a line bundle over the N\'eron 2-gon.
In this picture the two disjoint sections $C_0=C_0'\cup C_0''$,
$C_2=C_2'\cup C_2''$ are the zero section and the section at
infinity. Let $\foB\to\Delta$ be a smoothing of $B$. Since the line
bundle defining $Y$ extends to $\foB$ there exists an extension
$\foY\to \foB$ of $Y\to B$ with disjoint sections
$\foC_0,\foC_2\subset \foY$. The gluing $C_0\to C_2$ is given by an
automorphism of $B$ of order $2$. This automorphism can be extended
to an automorphism of $\foB/\Delta$ of the same order. An appropriate
identification of the sections now yields the desired deformation of
$X_2$.
\qed


\end{document}